\documentclass[a4paper,english]{article}

\usepackage[latin9]{inputenc}
\usepackage{geometry}
\usepackage{pgfplots}
\usepackage{float}
\usepackage{amsthm}
\usepackage{amsmath}
\usepackage{amstext}
\usepackage{amssymb}
\usepackage{color}
\usepackage{graphicx}
\usepackage{nicefrac}
\usepackage{fancyhdr}
\usepackage{lastpage}

\usepackage[english]{babel}
\usepackage{slashbox}
\usepackage{tikz}
\usepackage{graphicx,epstopdf}
\usepackage{xcolor}
\usepackage{rotfig}
\usepackage{rotating}
\usetikzlibrary{arrows,chains,matrix,positioning,scopes,shapes,arrows,calc,decorations.pathmorphing}
\usepackage{blkarray}
\usepackage{multirow}
 \usepackage{booktabs}
 \usepackage{algorithm}
\usepackage{hyperref}

\usepackage{algorithmic}
\usepackage{listings} 
\usepackage{tikz}
\usepackage{ifthen}
\usepackage{url}

\usepackage[font={small,it}]{caption}
\usepackage{subcaption}
\usepackage{amsthm}
\usepackage{thmtools}

\newcommand{\smb}{\left[\begin{smallmatrix}}
\newcommand{\sme}{\end{smallmatrix}\right]}

\newtheorem{theorem}{Theorem}

\newtheorem{example}[theorem]{Example}
 
\DeclareMathOperator{\h}{\mathcal{H}}

\DeclareMathOperator{\off}{off}

\newcommand{\R}{{\mathbb R}}

\newcommand{\Matlab}{{\sc Matlab}}

\newcommand{\calO}{{\mathcal O}}

\newcommand{\calH}{{\mathcal H}}

\newcommand{\calT}{{\mathcal T}}

\tikzstyle{block} = [rectangle, draw, fill=blue!20, 
    text width=5em, text centered, rounded corners, minimum height=4em]
    
\tikzset{
    ncbar angle/.initial=90,
    ncbar/.style={
        to path=(\tikztostart)
        -- ($(\tikztostart)!#1!\pgfkeysvalueof{/tikz/ncbar angle}:(\tikztotarget)$)
        -- ($(\tikztotarget)!($(\tikztostart)!#1!\pgfkeysvalueof{/tikz/ncbar angle}:(\tikztotarget)$)!\pgfkeysvalueof{/tikz/ncbar angle}:(\tikztostart)$)
        -- (\tikztotarget)
    },
    ncbar/.default=0.5cm,
}    

\definecolor{color1}{RGB}{137,207,240}

\tikzset{round left paren/.style={ncbar=0.5cm,out=100,in=-100}}
\tikzset{round right paren/.style={ncbar=0.5cm,out=80,in=-80}}

\title{Fast QR decomposition of HODLR matrices}
\author{Daniel Kressner\thanks{MATH-ANCHP, \'{E}cole Polytechnique F\'{e}d\'{e}rale de Lausanne, Station 8, 1015 Lausanne, Switzerland. E-mail: daniel.kressner@epfl.ch.} \and Ana 
\v{S}u\v{s}njara\thanks{MATH-ANCHP, \'{E}cole Polytechnique F\'{e}d\'{e}rale de Lausanne, Station 8, 1015 Lausanne, Switzerland. E-mail: susnjara.ana@gmail.com.  The work of Ana \v{S}u\v{s}njara
has been supported by the SNSF research project \emph{Low-rank updates of matrix functions and fast eigenvalue solvers.} }   }

\begin{document}

\date{}

\maketitle

\begin{abstract}
The efficient and accurate QR decomposition for matrices with hierarchical low-rank structures, such as HODLR and hierarchical matrices, has been challenging. Existing 
structure-exploiting algorithms are prone to numerical instability as they proceed indirectly, via Cholesky decompositions or a block Gram-Schmidt procedure. For a highly 
ill-conditioned matrix, such approaches either break down in finite-precision arithmetic or result in significant loss of orthogonality. Although these issues can sometimes be addressed by 
regularization and iterative refinement, it would be more desirable to have an algorithm that avoids these detours and is numerically robust to ill-conditioning. In this work, 
we propose such an algorithm for HODLR matrices. It achieves accuracy by utilizing Householder reflectors. It achieves efficiency by utilizing fast operations in 
the HODLR format in combination with compact WY representations and the recursive QR decomposition by Elmroth and Gustavson. Numerical experiments demonstrate that our 
newly proposed algorithm is robust to ill-conditioning and capable of achieving numerical orthogonality down to the level of roundoff error.

\end{abstract}

 
%
%
%
 \section{Introduction}

A HODLR (\emph{hierarchically off-diagonal low-rank}) matrix $A$ is defined recursively via $2\times 2$ block partitions of the form
\begin{equation} \label{eq:toplevelpartitioning}
 A =\left[ \begin{array}{c|c}
   A_{11} & A_{12} \\ \hline
   A_{21} & A_{22}\\ 
     \end{array} \right],
\end{equation}
where the off-diagonal blocks $A_{21}, A_{12}$ have low rank and the diagonal blocks are again HODLR matrices. The recursion is stopped once the diagonal blocks are of 
sufficiently small size, in the range of, say, a few hundreds. Storing the off-diagonal blocks in terms of their low-rank factors significantly reduces memory requirements and, 
potentially, the computational cost of operating with HODLR matrices. The goal of this work is to devise an efficient and numerically accurate algorithm for computing a QR decomposition
\[
 A = QR,
\]
where $R$ is an upper triangular HODLR matrix and the orthogonal matrix $Q$ is represented in terms of its so called compact WY 
representation~\cite{Schreiber1989}: $Q = I - YTY^T$, with the identity matrix $I$ and triangular/trapezoidal HODLR matrices $T,Y$.

HODLR matrices constitute one of the simplest data-sparse formats among the wide range of hierarchical low-rank formats that have been discussed in the literature during the last two decades. They have 
proved to be effective, for example, in solving large-scale linear systems~\cite{Aminfar2016} and operating with multivariate Gaussian distributions~\cite{Ambikasaran2016}. In our 
own work~\cite{Kressner2017,Susnjara2018} HODLR matrices have played a central role in developing fast algorithms for solving symmetric 
banded eigenvalue problems. In particular, a fast variant of the so called QDWH algorithm~\cite{NakaBaiGygi2010,NakaHigh2013} for computing spectral projectors requires the QR decomposition of a HODLR matrix. It is 
also useful for orthonormalizing data-sparse vectors. Possibly more importantly, the QR decomposition offers a stable alternative to the LU decomposition (without pivoting) for solving linear systems with 
nonsymmetric HODLR matrices or to the Cholesky decomposition applied to the normal equations for solving linear least-squares problems.

For the more general class of hierarchical matrices~\cite{Hackbusch2015}, a number of approaches aim at devising fast algorithms for QR decompositions~\cite{Bebendorf2008,BennerMach2010,Lintner2002}.
However, as we explain in Section~\ref{sec:overview} below, all existing approaches have limitations in terms of numerical accuracy and orthogonality, especially when $A$ is ill-conditioned. In 
the other direction, when further conditions are imposed on a HODLR matrix, leading to formats such as HSS (hierarchically semi-separable) or quasi-separable matrices, then it can be possible to 
devise QR or URV decompositions that fully preserve the structure; see~\cite{Eidelman2014,VandeVanBarelMastro2008,Xia2010} and the references therein. This is clearly not 
possible for HODLR matrices: In general, $Q$ and $R$ do not inherit from $A$ the property of having low-rank off-diagonal blocks. However, as it turns out, these factors 
can be very well \emph{approximated} via HODLR matrices. The approach presented in this work to obtain such approximations is different from any existing approach we are aware of. It is based on the recursive 
QR decomposition proposed by Elmroth and Gustavson~\cite{Elmroth2000} for dense matrices. The key insight in this work is that such a recursive algorithm
combines well with the use of the HODLR format for representing the involved compact WY representations. Demonstrated by the numerical experiments, the resulting algorithm is 
not only fast but it is also capable of yielding high accuracy, that is, a residual and orthogonality down to the level of roundoff error.

The rest of this paper is organized as follows. Section~\ref{sec:overview} provides an overview of HODLR matrices and the corresponding arithmetics, as well as the 
existing approaches to fast QR decompositions. In Section~\ref{sec:recursive_qr}, we recall the recursive QR decomposition from~\cite{Elmroth2000} for dense matrices. The 
central part of this work, Section~\ref{sec:hodlrqr} combines~\cite{Elmroth2000} with the HODLR format. Various numerical experiments 
reported in Section~\ref{sec:numericalresults} demonstrate the effectiveness of our approach. Finally, in Section~\ref{sec:rectangular}, we sketch the extension of our 
newly proposed approach from square to rectangular HODLR matrices.

\section{Overview of HODLR matrices and existing methods}
\label{sec:overview}

In the following, we focus our description on \emph{square} HODLR matrices. Although the extension of our algorithm to rectangular matrices does not require 
any substantially new ideas, the formal description of the algorithm would become significantly more technical. We have therefore postponed the rectangular case to Section~\ref{sec:rectangular}.

\subsection{HODLR matrices} \label{sec:hodlr}

\begin{figure}[ht!]
\centering
\begin{tikzpicture}[scale=0.75]
\fill[pink] (0, 0)--(1,0)--(1,8)--(0, 8);

\fill[gray] (0.5,7)--(1,7)--(1,7.5)--(0.5,7.5);
\fill[gray] (0,8)--(0,7.5)--(0.5,7.5)--(0.5,8);

\fill[gray] (1.5,6)--(2,6)--(2,6.5)--(1.5,6.5);
\fill[gray] (1,7)--(1,6.5)--(1.5,6.5)--(1.5,7);

\fill[gray] (2.5,5)--(3,5)--(3,5.5)--(2.5,5.5);
\fill[gray] (2,6)--(2,5.5)--(2.5,5.5)--(2.5,6);

\fill[gray] (3.5,4)--(4,4)--(4,4.5)--(3.5,4.5);
\fill[gray] (3,5)--(3,4.5)--(3.5,4.5)--(3.5,5);

\fill[gray] (4.5,3)--(5,3)--(5,3.5)--(4.5,3.5);
\fill[gray] (4,4)--(4,3.5)--(4.5,3.5)--(4.5,4);

\fill[gray] (5.5,2)--(6,2)--(6,2.5)--(5.5,2.5);
\fill[gray] (5,3)--(5,2.5)--(5.5,2.5)--(5.5,3);

\fill[gray] (6.5,1)--(7,1)--(7,1.5)--(6.5,1.5);
\fill[gray] (6,2)--(6,1.5)--(6.5,1.5)--(6.5,2);

\fill[gray] (7.5,0)--(8,0)--(8,0.5)--(7.5,0.5);
\fill[gray] (7,1)--(7,0.5)--(7.5,0.5)--(7.5,1);


\fill[gray] (0.1,0.1)--(0.1, 3.9)--(0.3,3.9)--(0.3, 0.1);
\fill[gray] (0.4,3.9)--(0.4, 3.7)--(3.9,3.7)--(3.9, 3.9);

\fill[gray] (0.1, 4.1)--(0.1,5.9)--(0.3, 5.9)--(0.3,4.1);
\fill[gray] (0.4, 5.9)--(0.4,5.7)--(1.9, 5.7)--(1.9,5.9);

\fill[gray] (0.1, 6.1)--(0.1,6.9)--(0.3, 6.9)--(0.3,6.1);
\fill[gray] (0.4, 6.9)--(0.4,6.7)--(0.9, 6.7)--(0.9,6.9);
 
\fill[gray] (0.1, 7.1)--(0.1,7.4)--(0.2, 7.4)--(0.2,7.1);
\fill[gray] (0.25, 7.4)--(0.25,7.3)--(0.45,7.3)--(0.45,7.4);

\fill[gray] (4.1,4.1)--(4.1, 7.9)--(4.3,7.9)--(4.3, 4.1);
\fill[gray] (4.4,7.9)--(4.4, 7.7)--(7.9,7.7)--(7.9, 7.9);

\fill[gray] (2.1, 6.1)--(2.1,7.9)--(2.3, 7.9)--(2.3,6.1);
\fill[gray] (2.4, 7.9)--(2.4,7.7)--(3.9, 7.7)--(3.9,7.9);

\fill[gray] (1.1, 7.1)--(1.1,7.9)--(1.3, 7.9)--(1.3,7.1);
\fill[gray] (1.4, 7.9)--(1.4,7.7)--(1.9, 7.7)--(1.9,7.9);

\fill[gray] (0.6, 7.6)--(0.6,7.9)--(0.7, 7.9)--(0.7,7.6);
\fill[gray] (0.75, 7.9)--(0.75,7.8)--(0.95,7.8)--(0.95,7.9);

\fill[gray] (1.6, 6.6)--(1.6,6.9)--(1.7, 6.9)--(1.7,6.6);
\fill[gray] (1.75, 6.9)--(1.75,6.8)--(1.95,6.8)--(1.95,6.9);

\fill[gray] (1.1, 6.1)--(1.1,6.4)--(1.2, 6.4)--(1.2,6.1);
\fill[gray] (1.25, 6.4)--(1.25,6.3)--(1.45,6.3)--(1.45,6.4);

\fill[gray] (2.1, 5.1)--(2.1,5.4)--(2.2, 5.4)--(2.2,5.1);
\fill[gray] (2.25, 5.4)--(2.25,5.3)--(2.45,5.3)--(2.45,5.4);

\fill[gray] (3.1, 4.1)--(3.1,4.4)--(3.2, 4.4)--(3.2,4.1);
\fill[gray] (3.25, 4.4)--(3.25,4.3)--(3.45,4.3)--(3.45,4.4);

\fill[gray] (4.1, 3.1)--(4.1,3.4)--(4.2, 3.4)--(4.2,3.1);
\fill[gray] (4.25, 3.4)--(4.25,3.3)--(4.45,3.3)--(4.45,3.4);

\fill[gray] (5.1, 2.1)--(5.1,2.4)--(5.2, 2.4)--(5.2,2.1);
\fill[gray] (5.25, 2.4)--(5.25,2.3)--(5.45,2.3)--(5.45,2.4);

\fill[gray] (6.1, 1.1)--(6.1,1.4)--(6.2, 1.4)--(6.2,1.1);
\fill[gray] (6.25, 1.4)--(6.25,1.3)--(6.45,1.3)--(6.45,1.4);

\fill[gray] (7.1, 0.1)--(7.1,0.4)--(7.2, 0.4)--(7.2,0.1);
\fill[gray] (7.25, 0.4)--(7.25,0.3)--(7.45,0.3)--(7.45,0.4);

\fill[gray] (2.6, 5.6)--(2.6,5.9)--(2.7, 5.9)--(2.7,5.6);
\fill[gray] (2.75, 5.9)--(2.75,5.8)--(2.95,5.8)--(2.95,5.9);

\fill[gray] (3.6, 4.6)--(3.6,4.9)--(3.7, 4.9)--(3.7,4.6);
\fill[gray] (3.75, 4.9)--(3.75,4.8)--(3.95,4.8)--(3.95,4.9);

\fill[gray] (4.6, 3.6)--(4.6,3.9)--(4.7, 3.9)--(4.7,3.6);
\fill[gray] (4.75, 3.9)--(4.75,3.8)--(4.95,3.8)--(4.95,3.9);

\fill[gray] (5.6, 2.6)--(5.6,2.9)--(5.7, 2.9)--(5.7,2.6);
\fill[gray] (5.75, 2.9)--(5.75,2.8)--(5.95,2.8)--(5.95,2.9);

\fill[gray] (6.6, 1.6)--(6.6,1.9)--(6.7, 1.9)--(6.7,1.6);
\fill[gray] (6.75, 1.9)--(6.75,1.8)--(6.95,1.8)--(6.95,1.9);

\fill[gray] (7.6, 0.6)--(7.6,0.9)--(7.7, 0.9)--(7.7,0.6);
\fill[gray] (7.75, 0.9)--(7.75,0.8)--(7.95,0.8)--(7.95,0.9);

\fill[gray] (4.1, 0.1)--(4.1,1.9)--(4.3, 1.9)--(4.3,0.1);
\fill[gray] (4.4, 1.9)--(4.4,1.7)--(5.9, 1.7)--(5.9,1.9);

\fill[gray] (4.1, 0.1)--(4.1,1.9)--(4.3, 1.9)--(4.3,0.1);
\fill[gray] (4.4, 1.9)--(4.4,1.7)--(5.9, 1.7)--(5.9,1.9);

\fill[gray] (6.1, 2.1)--(6.1,3.9)--(6.3, 3.9)--(6.3,2.1);
\fill[gray] (6.4, 3.9)--(6.4,3.7)--(7.9, 3.7)--(7.9,3.9);


\fill[gray] (2.1, 4.1)--(2.1,4.9)--(2.3, 4.9)--(2.3,4.1);
\fill[gray] (2.4, 4.9)--(2.4,4.7)--(2.9, 4.7)--(2.9,4.9);

\fill[gray] (4.1, 2.1)--(4.1,2.9)--(4.3, 2.9)--(4.3,2.1);
\fill[gray] (4.4, 2.9)--(4.4,2.7)--(4.9, 2.7)--(4.9,2.9);

\fill[gray] (6.1, 0.1)--(6.1, 0.9)--(6.3, 0.9)--(6.3,0.1);
\fill[gray] (6.4, 0.9)--(6.4,0.7)--(6.9, 0.7)--(6.9,0.9);


\fill[gray] (3.1, 5.1)--(3.1,5.9)--(3.3, 5.9)--(3.3,5.1);
\fill[gray] (3.4, 5.9)--(3.4,5.7)--(3.9, 5.7)--(3.9,5.9);

\fill[gray] (5.1, 3.1)--(5.1,3.9)--(5.3, 3.9)--(5.3,3.1);
\fill[gray] (5.4, 3.9)--(5.4,3.7)--(5.9, 3.7)--(5.9,3.9);

\fill[gray] (7.1, 1.1)--(7.1,1.9)--(7.3, 1.9)--(7.3,1.1);
\fill[gray] (7.4, 1.9)--(7.4,1.7)--(7.9, 1.7)--(7.9,1.9);


\draw (0, 0)--(1,0)--(1,8)--(0, 8)--(0,0);
\draw (0,0) rectangle (8,8);
\draw (0,0) rectangle  (4, 4);
\draw (4,0) rectangle  (8, 4);
\draw (0,4) rectangle  (4, 8);

\draw (0,6)--(4,6);
\draw (2,4)--(2,8);

\draw (4,2)--(8,2);
\draw (6,0)--(6,4);

\draw (0,7)--(2,7);
\draw (1,6)--(1,8);

\draw (0,7.5)--(1,7.5);
\draw (0.5,7)--(0.5,8);

\draw (1,6.5)--(2,6.5);
\draw (1.5,6)--(1.5,7);

\draw (3,4)--(3,6);
\draw (2,5)--(4,5);

\draw (2,5.5)--(3,5.5);
\draw(2.5,5)--(2.5,6);

\draw (3,4.5)--(4,4.5);
\draw (3.5,4)--(3.5,5);

\draw (4,3)--(6,3);
\draw (5,2)--(5,4);

\draw (4,3.5)--(5,3.5);
\draw (4.5,3)--(4.5,4);

\draw[] (5,2.5)--(6,2.5);
\draw (5.5,2)--(5.5,3);

\draw (6,1)--(8,1);
\draw (7,0)--(7,2);

\draw (6,1.5)--(7,1.5);
\draw (6.5,1)--(6.5,2);

\draw (7,0.5)--(8,0.5);
\draw (7.5,0)--(7.5,1);
\end{tikzpicture} 

\caption{HODLR matrix of level $\ell = 4$. The high-lighted block column is processed recursively by our newly proposed algorithm; see Section~\ref{sec:hodlrqr}. It consists of 
three different types of blocks: On top a level-one HODLR matrix, below a low-rank matrix entirely contained within the high-lighted block column, and two low-rank matrices extending into other block columns. }
\label{fig:hodlr_matrix}
\end{figure}
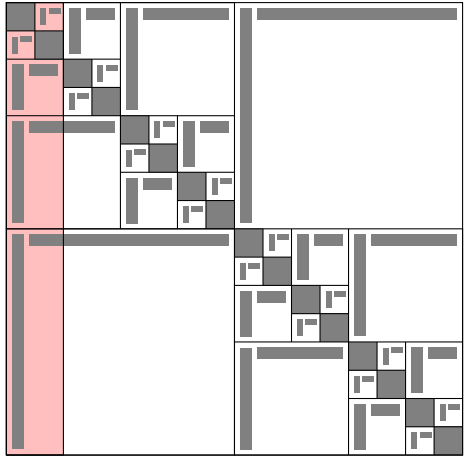

As discussed in the introduction, a HODLR matrix $A\in \R^{n\times n}$ is defined by performing a recursive partition of the form~\eqref{eq:toplevelpartitioning} and requiring 
all occuring off-diagonal blocks to be of low rank. When this recursion is performed $\ell$ times, we say that $A$ is a HODLR matrix of \emph{level $\ell$}; see Figure~\ref{fig:hodlr_matrix} for an illustration.

Clearly, the definition of a HODLR matrix depends on the block sizes chosen in~\eqref{eq:toplevelpartitioning} on every level of the recursion or, equivalently, on the integer partition
\begin{equation} \label{eq:integerpart}
 n = n_1 + n_2 + \cdots n_{2^\ell},  
\end{equation}
defined by the sizes $n_j\times n_j$, $j = 1,\ldots, 2^\ell,$ of the diagonal blocks on the lowest level of the recursion. If possible, it is advisable to choose the 
level $\ell$ and the integers $n_j$ such that all $n_j$ are nearly equal to a prescribed minimal block size $n_{\min}$. In the following, when discussing the complexity of operations, we assume that such a balanced partition has been chosen.

Given an integer partition~\eqref{eq:integerpart}, we define $\mathcal{H}_{n\times n}(\ell, k)$ to be the set of $n \times n$ HODLR  matrices of level $\ell$ and rank (at most) $k$, that is, $A \in \mathcal{H}_{n\times n}(\ell, k)$ if every off-diagonal block of $A$ in 
the recursive block partition induced by~\eqref{eq:integerpart} has rank at most $k$.

A matrix $A \in \calH_{n\times n}(\ell,k)$ admits a data-sparse representation by storing its off-diagonal blocks in 
terms of their low-rank factors. Specifically, letting $A\vert_{\off} \in \R^{n_L\times n_R}$ denote an arbitrary off-diagonal block in the recursive partition of $A$, we can write
\begin{equation} \label{eq:factoroffdiagonal}
 A\vert_{\off} =A_L A_R, \quad A_L \in \R^{n_L \times k}, \quad A_R \in \R^{k\times n_R}.
\end{equation}
Storing $A_L$ and $A_R$ instead of $A\vert_{\off}$ for every such off-diagonal block reduces the overall memory required for storing $A$ from $\mathcal O(n^2)$ to 
$\mathcal{O}(kn\log n)$.

We call a factorization~\eqref{eq:factoroffdiagonal} \emph{left-orthogonal} if $A_L^T A_L = I_k$. Provided that $k \le n_L$, an arbitrary factorization~\eqref{eq:factoroffdiagonal} can be turned
into a left-orthogonal one by computing an economy-sized QR decomposition $A_L = QR$, see~\cite[Theorem 5.2.3]{GolVanL2013}, and replacing
$A_L \gets Q$, $A_R \gets R A_R$.  This described procedure requires $\calO((n_L + n_R)k^2)$ operations.   

\subsubsection{Approximation by HODLR matrices}
\label{sec:approximation_HODLR}

A general matrix $A \in \mathbb R^{n\times n}$ can be approximated by a HODLR matrix by performing low-rank truncations of the off-diagonal blocks in the recursive block partition. Specifically, 
letting $A\vert_{\off}$ denote such an off-diagonal block, one computes a singular value decomposition $A\vert_{\off} = U\Sigma V^T$ with the diagonal 
matrix $\Sigma = \text{diag}(\sigma_1,\sigma_2, \ldots)$ containing the singular values. Letting $U_k$ and $V_k$ contain the first $k$ columns of $U$ and $V$, respectively, and 
setting $\Sigma_k = \text{diag}(\sigma_1,\ldots,\sigma_k)$, one obtains 
a rank-$k$ approximation
\begin{equation}\label{eq:lrapprox}
A\vert_{\off} \approx A_L A_R
\end{equation}
by setting $A_L = U_k$ and $A_R = \Sigma_kV_k^T$. We note that this approximation is optimal among all rank-$k$ matrices for any unitarily 
invariant norm~\cite[Section 7.4.9]{Horn2013}. In particular, for the matrix $2$-norm, we have
\begin{equation} \label{eq:localapproximationerror}
 \|A\vert_{\off} - A_L A_R\|_2 = \sigma_{k+1}.
\end{equation}
In passing, we note that the factorization chosen in~\eqref{eq:lrapprox} is left-orthogonal.

Performing the approximation~\eqref{eq:lrapprox} for every off-diagonal block in the recursive block partition yields a HODLR matrix $A_{\h,k} \in \mathcal{H}_{n\times n}(\ell, k)$. 

In practice, the rank $k$ is chosen adaptively and separately for each off-diagonal block $A\vert_{\off}$. Given a prescribed tolerance $\epsilon > 0$, we 
choose $k \equiv k_\epsilon$ to be the smallest integer such that $\sigma_{k_\epsilon+1} \le \epsilon$. In turn,~\eqref{eq:localapproximationerror} implies
$\|A\vert_{\off} - A_L A_R\|_2 \le \epsilon$. The resulting HODLR approximation $A_{\h, \epsilon}$ satisfies $\Vert A - A_{\h, \epsilon}\Vert_2 \leq \ell \epsilon$; see, e.g.,~\cite[Theorem 2.2]{Bini2017}.

\paragraph{Recompression.} Most manipulations involving HODLR matrices lead to an increase of off-diagonal ranks. This increase is potentially mitigated by performing recompression.
Let us consider an off-diagonal block $A\vert_{\off}= A_L A_R$, with $A_L \in \R^{n_L\times k_A}, A_R \in \R^{k_A\times n_R}$ and choose $k_\epsilon$ as explained above. If $k_\epsilon < k_A$, a 
rank-$k_\epsilon$ approximation reduces memory requirements while maintaining $\epsilon$-accuracy. We use the following well-known procedure for effecting this approximation.
\begin{enumerate}
 \item Compute $A_L = Q_1 R_1$ and $A_R^T = Q_2 R_2$, economy-sized QR decompositions of $A_L$ and $A^T_R$, respectively.
 \item Compute SVD $R_1R_2^T = \tilde{U} \Sigma \tilde{V}^T$.
 \item Update $A_L \gets Q_1\tilde{U}_{k_\epsilon}$ and $A_R \gets \Sigma_{k_\epsilon} \tilde{V}_{k_\epsilon}^T Q_2^T$.
\end{enumerate}
This procedure, which will be denoted by $\calT_\epsilon$, requires $\calO\left((n_L+n_R)k_A^2\right)$ operations.

\subsubsection{Operating with HODLR matrices} \label{sec:HODLR}

A number of operations can be performed efficiently with HODLR matrices. Table~\ref{table:complexity_HODLR} lists the operations relevant in this work, together with their computational complexity; see, e.g.,~\cite[Chapter 3]{Hackbusch2015} for 
more details. It is important to note that all operations, except for matrix-vector multiplication, are combined 
with low-rank truncation, as discussed above, to limit rank growth in the off-diagonal blocks. The symbol $\h$ signifies the inexactness due to truncations. The 
complexity estimates assume that all off-diagonal ranks encountered during an operation remain $\calO(k)$.
 
\begin{table}[ht!]
 \caption{Complexity of operations with HODLR matrices: $A_1, \ldots, A_6 \in \mathcal H_{n\times n}(\ell,k)$,  with $A_3$ invertible,  $A_4$ invertible upper triangular, $A_6$ symmetric positive definite,   
 $ U,V\in \R^{n\times p}$ with $p = \mathcal O(k)$, and $v\in \R^{n}$. }
\label{table:complexity_HODLR}
{\renewcommand{\arraystretch}{1.1}
\centering \begin{tabular}{rcl}
 \hline
 Operation &  &Computational complexity \\ \hline
 Matrix-vector multiplication:   $A_1 v$ &  & $\mathcal{O}(kn\log n)$ \\
 Matrix addition: $A_1 +_{\mathcal{H}} A_2$  &  & $\mathcal{O}(k^2n \log n)$ \\
 Matrix low-rank update: $A_1 +_{\h} UV^T$  &  & $\mathcal{O}(k^2 n\log n)$  \\
 Matrix-matrix multiplication: $A_1 *_{\mathcal{H}} A_2$ &  & $\mathcal{O}(k^2 n\log^2 n)$  \\
 Matrix inversion: $\h\operatorname{-inv}(A_3)$ & & $\mathcal{O}(k^2 n\log^2 n)$ \\
 Solution of triangular matrix equation: $A_5 *_{\mathcal{H}} A^{-1}_4$ & &$\mathcal{O}(k^2 n\log^2 n)$  \\ 
 Cholesky decomposition: $\h\operatorname{-Cholesky}(A_6)$ & & $\mathcal{O}(k^2 n\log^2 n)$  \\ \hline
 \end{tabular}
 }
 \end{table}

\subsection{Cholesky-based QR decomposition} \label{sec:cholesky}

This and the following sections describe three existing methods for efficiently computing the QR decomposition of a HODLR matrix. All these methods have originally been proposed 
for the broader class of hierarchical matrices.

The first method, proposed by Lintner~\cite{Lintner2002,Lintner2004}, is based on the well-known connection between the QR and Cholesky decompositions. Specifically, letting $A = QR$ 
be the QR decomposition of an invertible $n\times n$ matrix $A$, we have 
\[ A^TA = R^T Q^TQR = R^TR. \]
Thus, the upper triangular factor $R$ can be obtained from the Cholesky decomposition of the symmetric positive definite matrix $A^T A$.
The orthogonal factor $Q$ is obtained from solving the triangular system $A = QR$. 
According to Table~\ref{table:complexity_HODLR}, these three steps (forming $A^T A$, computing the Cholesky decomposition, solving the triangular matrix equation) require 
$\calO(k^2n \log^2n)$ operations in the HODLR format.

For dense matrices, the approach described above is well-known and often called \emph{CholeskyQR algorithm}; see~\cite[Pg. 214]{Stewart1973} for an early reference. A major disadvantage of this approach, $Q$ rapidly loses orthogonality in finite precision arithmetic as the condition number of $A$ increases.
As noted in~\cite{Stathopoulos2002}, the numerical orthogonality $\|Q^T Q - I\|_2$ is usually at the level of the squared condition 
number $\kappa(A^T A) = \kappa(A)^2$ times the unit roundoff $\mathsf u$. To improve its orthogonality, one can apply the CholeskyQR algorithm again to $Q$ and 
update $R$ accordingly. As shown in~\cite{Yamamoto2015}, this so called \emph{CholeskyQR2 algorithm} results in a numerically orthogonal factor, provided that $\kappa(A)$ is at most $\calO(\mathsf u^{-1/2})$.

The CholeskyQR2 algorithm for HODLR and hierarchical matrices~\cite{Lintner2002} is additionally affected by low-rank truncation and may require several reorthogonalization steps to reach numerical orthogonality on the level of the truncation error, increasing the computational cost. Another approach proposed in~\cite{Lintner2002} to avoid loss of 
orthogonality is to first compute a  polar decomposition $A = QH$ and then apply the CholeskyQR algorithm to $H$. Because of 
$\kappa(H) = \kappa(A) = \sqrt{\kappa(A^T A)}$, this improves the accuracy of the CholeskyQR algorithm. On the other hand, the need for computing the polar 
decomposition via an iterative method, such as the sign-function iteration~\cite{Higham1986}, also significantly increases the computational cost.

\subsection{LU-based QR decomposition}

An approach proposed by Bebendorf~\cite[Sec. 2.10]{Bebendorf2008} can be viewed as orthogonalizing a recursive block LU decomposition. Given $A\in \R^{n \times n}$, let us  
partition \begin{equation} \label{eq:partition2}
 A =\left[ \begin{array}{c|c}
   A_{11} & A_{12} \\ \hline
   A_{21} & A_{22}\\ 
     \end{array} \right]
\end{equation}
and suppose that $A_{11}$ is invertible. Setting $X = A_{21}A^{-1}_{11}$, consider the block LU decomposition
\begin{equation} \label{eq:blocklu}
 A = \begin{bmatrix}
   I & 0 \\ 
   X & I \\ 
     \end{bmatrix}\begin{bmatrix}
   A_{11} & A_{12} \\ 
   0 & A_{22} - X A_{12} \\ 
     \end{bmatrix}.
\end{equation}
The first factor is orthogonalized by (1) rescaling the first block column with the inverted Cholesky factor of the symmetric positive definite matrix
$I + X^TX = R_1^T R_1$ and (2) choosing the second block column as $\big[ { -X^T \atop I} \big]$, scaled with 
the inverted Cholesky factor of $I + XX^T = R_2^T R_2$.
Adjusting the second factor in~\eqref{eq:blocklu} accordingly does not change its block triangular structure. More precisely, one can prove that
\[
 A = \underbrace{\begin{bmatrix} I & -X^T\\ X& I\end{bmatrix} \begin{bmatrix} R^{-1}_1  &0\\ 0& R^{-1}_{2}\end{bmatrix}}_{=:\tilde Q}
 \underbrace{\begin{bmatrix} R_1A_{11}  &R^{-T}_1(A_{12} + X^TA_{22})\\ 0& R^{-T}_{2}(A_{22} - XA_{12})\end{bmatrix}}_{=:\tilde R}
\]
holds. By construction, $\tilde Q$ is orthogonal. This procedure is applied recursively to the diagonal blocks of $\tilde R$. If $A$ is a HODLR matrix corresponding to the partition~\eqref{eq:partition2} at 
every level of the recursion then all involved operations can be performed efficiently in the HODLR format. Following~\cite{Bebendorf2008}, the overall computational cost is, once again, $\calO(k^2n\log^2 n)$. 

An obvious disadvantage of the described approach, it requires the leading diagonal block $A_{11}$ to be well conditioned for every subproblem encountered during the recursion. This rather 
restrictive assumption is only guaranteed for specific matrix classes, such as well-conditioned positive definite matrices.

\subsection{QR decomposition based on a block Gram-Schmidt procedure}

The equivalence between the QR decomposition and the Gram-Schmidt procedure for full-rank matrices is well known. In particular, applying the modified block Gram-Schmidt procedure to 
the columns of $A$ leads to the block recursive QR decomposition presented in~\cite[Sec. 5.2.4]{GolVanL2013}. Benner and Mach~\cite{BennerMach2010} combined this idea with hierarchical 
matrix arithmetic. In the following, we briefly summarize their approach. Partitioning the economy-sized QR decomposition of $A$ into block columns yields the relation
\begin{equation} \label{eq:blockcolpartition}
[A_1 \hskip 3pt A_2] = [Q_1 \hskip 3pt Q_2] \begin{bmatrix} R_{11} &R_{12} \\ 0 &R_{22} \end{bmatrix}. 
\end{equation}
This 
yields
three steps for the computation of $Q$ and $R$:
\begin{enumerate}
 \item Compute (recursively) the QR decomposition $A_1 = Q_1R_{11}$.
\item  Compute $R_{12} = Q_1^TA_2$ and update $A_2 \gets A_2 - Q_1R_{12}$.
\item Compute (recursively) the QR decomposition  $A_2 = Q_2R_{22}$.
\end{enumerate}
Step 2 can be implemented efficiently for a HODLR matrix that aligns with the block column partitioning~\eqref{eq:blockcolpartition}. Steps 1 and 3 are executed 
recursively until the lowest level of the HODLR structure is reached. On this lowest level, it is suggested in~\cite[Alg. 3]{BennerMach2010} to compute the QR decomposition of 
compressed block columns. We refrain from providing details and point out that we consider similarly compressed block columns in Section~\ref{sec:hodlrqr} below. The overall 
computational complexity is $\calO(k^2n\log^2 n)$.

The described algorithm inherits the numerical instability of Gram-Schmidt procedures. In particular, we cannot expect to obtain a
numerically orthogonal factor $Q$ in finite-precision arithmetic when $A$ is ill-conditioned; see also the analysis in~\cite[Sec. 3.4]{BennerMach2010}.

\section{Recursive WY-based QR decomposition}
\label{sec:recursive_qr}

In this section, we recall the recursive QR decomposition by Elmroth and Gustavson~\cite{Elmroth2000} for a general, dense $m\times n$ matrix $A$ with $m\ge n$. The orthogonal factor $Q$ is returned in terms of the compact WY representation~\cite{Schreiber1989} of the $n$ Householder reflectors involved in the decomposition:
\begin{equation} \label{eq:compactwy}
 Q = I_m - Y T Y^T,
\end{equation}

where $T$ is an $n\times n$ upper triangular matrix and $Y$ is an $m\times n$ matrix with the first $n$ rows in unit lower triangular form.

For $n = 1$, the matrix $A$ becomes a column vector and we let $Q = I_m - \gamma y y^T$ be the Householder reflector~\cite[Sec. 5.1.2]{GolVanL2013} that maps $A$ to a scalar multiple of the unit vector. Then $Q$ is trivially of the form~\eqref{eq:compactwy}.

For $n > 1$, we partition $A$ into two block columns of roughly equal size:
\[
 A = \left[ \begin{array}{c|c}
      A_1 & A_2
     \end{array}\right], \quad A_1\in \R^{m\times n_1}, \quad A_2\in \R^{m\times n_2},\quad n = n_1 + n_2.
\]
By recursion, we compute a QR decomposition of the first block column
\[
 A_1 = Q_1 \begin{bmatrix} R_1 \\ 0 \end{bmatrix},\quad Q_1 = I_m - Y_1 T_1 Y_1^T,
\]
with $T_1\in \R^{n_1\times n_1}$, $Y_1\in \R^{m\times n_1}$ taking the form explained above. The second block column $A_2$ is updated,
\[
 \tilde A_2 = Q_1^T A_2 = A_2 - Y_1 T_1 (Y_1^T A_2),
\]
and then partitioned as
\[
 \tilde A_2 = \begin{bmatrix}
               A_{12} \\
               A_{22}
              \end{bmatrix}, \quad A_{12}\in \R^{n_1\times n_2}, \quad A_{22}\in \R^{(m-n_1)\times n_2}.
\]
Again by recursion, we compute a QR decomposition of the bottom block:
\[
 A_{22} = Q_2 \begin{bmatrix} R_2 \\ 0 \end{bmatrix},\quad Q_2 = I_{m-n_1} - Y_2 T_2 Y_2^T.
\]

To combine the QR decompositions of the first and the updated second block column, we embedd $Q_2$ into the larger matrix 
\[
 \tilde Q_2 = \begin{bmatrix} I_{n_1} & 0 \\
           0 & Q_2
          \end{bmatrix} = I_m - \tilde Y_2 T_2 \tilde Y_2^T, \quad \tilde Y_2 = \begin{bmatrix} 0 \\ Y_2 \end{bmatrix}.
\]
By setting
\[
 R = \left[ \begin{array}{c|c}
      R_{1} & A_{12} \\ \hline
      0 & R_2 \\
      \end{array} \right]
\]
and \begin{eqnarray}
 Q&=& Q_1 \tilde Q_2
= \big(I_m - Y_1 T_1 Y_1^T\big)\big( I_m - \tilde Y_2 T_2 \tilde Y_2^T  \big)  \nonumber \\
&=& I_m - Y_1 T_1 Y_1^T - \tilde Y_2 T_2\tilde Y_2^T + Y_1 T_1 Y_1^T \tilde Y_2 T_2 \tilde Y_2^T  \nonumber \\
&=& I_m - \left[ \begin{array}{c|c} Y_1 & \tilde Y_2 \end{array}\right]
\left[ \begin{array}{c|c} T_1 & -T_1 Y_1^T \tilde Y_2 T_2 \\ \hline 0 & T_2 \end{array}\right]
\left[ \begin{array}{c|c} Y_1 & \tilde Y_2 \end{array}\right]^T, \label{eq:recqrcombined}
\end{eqnarray}
we obtain a QR decomposition
\[
 A = Q \begin{bmatrix} R \\ 0 \end{bmatrix},\quad Q = I-YTY^T,\quad Y = \left[ \begin{array}{c|c} Y_1 & \tilde Y_2 \end{array}\right], \quad 
 T = \left[ \begin{array}{c|c} T_1 & -T_1 Y_1^T \tilde Y_2 T_2 \\ \hline 0 & T_2 \end{array}\right].
\]

Algorithm~\ref{alg:block_denseQR} summarizes the described procedure. To simplify the description, the recursion is performed down to individual columns. In practice~\cite{Elmroth2000}, the recursion is stopped earlier: When the number of columns does not exceed a certain block size $n_b$ (e.g., $n_b = 32$),
a standard Householder-based QR decomposition is used.

\begin{algorithm}[h!]
    \caption{\text{Recursive block QR decomposition}}
    \label{alg:block_denseQR}
    \textbf{Input:} Matrix $A \in \R^{m\times n}$ with $m\geq n$.\\
    \textbf{Output:} Matrices $Y \in \R^{m\times n}$, $T\in \R^{n\times n}$, $R\in \R^{n\times n}$, defining a QR decomposition $A = Q\begin{bmatrix}R\\ 0\end{bmatrix}$ with $Q = I_m - YTY^T \in \R^{m\times m}$ orthogonal. 

\begin{algorithmic}[1]
\STATE \textbf{function} $[Y, T, R] = {\tt{blockQR}} (A)$

\IF {$n = 1$ }
\STATE Compute Householder reflector $I_m - \gamma yy^T$ such that $(I_m - \gamma yy^T) A = \begin{bmatrix}\rho \\ 0\end{bmatrix}$.
\STATE Set $Y = y$, $T = \gamma$ and $R = \rho$.
\ELSE 
\STATE Set $n_1 = \lfloor n/2 \rfloor$. 
\STATE Call $[Y_1, T_{1}, R_{1}] = {\tt blockQR}(A(:, 1:n_1))$. 
\STATE Update $A(:, n_1+1:n) \gets (I - Y_1TY_1^T)^T A(:, n_1+1:n)$. 
\STATE Set $[Y_2, T_{2}, R_{2}] = {\tt blockQR}(A(n_1+1:m, n_1+1:n))$. 
\STATE Set $\tilde{Y}_2 =  \begin{bmatrix}0 \\ Y_2\end{bmatrix}$ and compute $T_{12} =-T_{1} Y_1^T \tilde{Y}_2 T_{2}$. 
\STATE Return  $Y = \begin{bmatrix} Y_1 & \tilde{Y}_2\end{bmatrix}$, $T = \begin{bmatrix}T_{1} &T_{12} \\ 0 & T_{2} \end{bmatrix}$ and
$R = \begin{bmatrix}R_{1} &A(1:n_1, n_1+1:n) \\ 0 & R_{2} \end{bmatrix}$. 
\ENDIF
\STATE \textbf{end function}
\end{algorithmic}
\end{algorithm}

\section{Recursive WY-based QR decomposition of HODLR matrices} \label{sec:hodlrqr}

By combining the recursive block QR decomposition (Algorithm~\ref{alg:block_denseQR}) with HODLR arithmetic, we will show in this section how to derive an efficient algorithm for computing the QR decomposition of a level-$\ell$ HODLR matrix $A \in \R^{n\times n}$.

The matrix processed in one step of the recursion of our algorithm takes the following form:
\begin{equation} \label{eq:defH}
H = 
\begin{bmatrix}
 \tilde A \\
 B \\
 C \\
\end{bmatrix},
\end{equation}
where:
\begin{itemize}
 \item $\tilde A\in \R^{m \times m}$ is a HODLR matrix of level $\tilde \ell \le \ell$;
 \item $B \in \R^{p \times m}$ is given in factorized form $B = B_L B_R$ with $B_L \in \R^{p \times r_1}$ and $B_R \in \R^{r_1 \times m}$ for some (small) integer $r_1$;
 \item $C\in \R^{r_2 \times m}$ for some (small) integer $r_2$.
\end{itemize}
To motivate this structure, it is helpful to consider the block column highlighted in Figure~\ref{fig:hodlr_matrix}. The first block in this block column consists of a (square) level-one HODLR matrix, corresponding 
to the matrix $\tilde A$ in~\eqref{eq:defH}. The other three blocks are all of low rank, and the matrix $B$ in~\eqref{eq:defH} corresponds to the first of these blocks. The bottom two blocks extend into the next
block column(s). For these two blocks, it is assumed that their left factors are orthonormal. These left factors are ignored and the parts of the right factors residing in the highlighted block column are collected in the matrix $C$ in~\eqref{eq:defH}.

Given a matrix $H$ of the form~\eqref{eq:defH}, we aim at computing, recursively and approximately, a QR decomposition of the form
\begin{equation} \label{eq:qrdecomp}
 H = Q \begin{bmatrix} R \\ 0 \end{bmatrix} , \quad Q = I - YTY^T
\end{equation}
such that $R,T \in \R^{m \times m}$ are upper triangular HODLR matrices of level $\tilde \ell$ and the structure of $Y$ reflects the structure of $H$, that is,
\begin{equation} \label{eq:defY}
Y = 
\begin{bmatrix}
 Y_A \\
 Y_B \\
 Y_C \\
\end{bmatrix},
\end{equation}
where: $Y_A \in \R^{m\times m}$ is a unit lower triangular HODLR matrix of level $\tilde \ell$; 
$Y_B \in \R^{p\times m}$ is in factorized form; and
$Y_C\in \R^{r_2 \times m}$.

On the highest level of the recursion, when $\tilde \ell = \ell$, the matrices $B,C$ in~\eqref{eq:defH} vanish, and $H = A$ is the original HODLR matrix we aim at decomposing. The QR decomposition returned
on the highest level has the form~\eqref{eq:qrdecomp} with both $Y,T$ triangular level-$\ell$ HODLR matrices.

The computation of the QR decomposition~\eqref{eq:qrdecomp} proceeds in several steps, which are detailed in the following.
\begin{paragraph}{Preprocessing.}
Using the procedure described in Section~\ref{sec:HODLR}, we may assume that the factorization of $B$ is normalized such that $B_L$ has orthonormal columns. 
For the moment, we will discard $B_L$ and aim at decomposing instead of $H$ the compressed matrix
\begin{equation} \label{eq:compressedcol}
\tilde H = \begin{bmatrix}
 \tilde A \\
 B_R \\
 C \\
\end{bmatrix},
\end{equation}
which has size $(m+r_1+r_2)\times m$.
\end{paragraph}

\begin{paragraph}{QR decomposition of $\tilde H$ on the lowest level, $\tilde \ell = 0$.}
On the lowest level of recursion, $\tilde A$ becomes a dense matrix. We perform a dense QR decomposition of the matrix $\tilde H$ defined in~\eqref{eq:compressedcol}. For this purpose, one can use, for example, Algorithm~\ref{alg:block_denseQR}. This yields the orthogonal factor $\tilde Q$ in terms of its compact WY representation, which we partition as
\begin{equation}
\label{eq:Q_tilde}
 \tilde Q = I - \tilde Y T \tilde Y^T, \qquad \tilde Y = \begin{bmatrix}
             Y_A \\
             \tilde Y_B \\
             Y_C
\end{bmatrix}, \quad Y_A \in \R^{m\times m}, \quad \tilde Y_B\in \R^{r_1\times m}, \quad Y_C\in \R^{r_2\times m}\text{.}
\end{equation}
\end{paragraph}

\begin{paragraph}{QR decomposition of $\tilde H$ on higher levels, $\tilde \ell \ge 1$.}
We proceed recursively as follows. First, $\tilde H$ is repartitioned as follows:
\begin{equation} \label{eq:repartH}
 \tilde H = \begin{bmatrix}
 \tilde A_{11} & \tilde A_{12} \\
 \tilde A_{21} & \tilde A_{22} \\
  B_{R,1} &  B_{R,2} \\
  C_1 & C_2
\end{bmatrix}.
\end{equation}
Here, $\tilde A = \begin{bmatrix}
\tilde A_{11} & \tilde A_{12} \\
 \tilde A_{21} & \tilde A_{22}
                  \end{bmatrix}$
is split according to its HODLR format, that is,  $\tilde A_{11}\in \R^{m_1\times m_1}$, $\tilde A_{22}\in \R^{m_2\times m_2}$, with $m = m_1 + m_2$, are 
HODLR matrices of level $\tilde{\ell}-1$, and $\tilde A_{21}, \tilde A_{12}$ are low-rank matrices stored in factorized form. 

Note that the first block column of $\tilde H$ in~\eqref{eq:repartH} has precisely the form~\eqref{eq:defH} with the level of the HODLR matrix reduced by one, the 
low-rank block given by $\tilde A_{21}$ and the dense part given by $\begin{bmatrix} B_{R,1} \\ C_1 \end{bmatrix}$. This allows us to apply recursion and obtain a QR decomposition
\begin{equation}\label{eq:first_block_QR}
 \begin{bmatrix}
 \tilde A_{11}  \\
 \tilde A_{21}  \\
 B_{R,1} \\
  C_1 \\
\end{bmatrix} = Q_1 \begin{bmatrix}
  R_1 \\
  0 \\ 
\end{bmatrix}, \quad Q_1 = I - Y_1 T_1 Y_1^T, \quad Y_1 = \begin{bmatrix} Y_{A,11} \\ Y_{A,21} \\ Y_{B_R,1} \\ Y_{C,1} \end{bmatrix},
\end{equation}
with a HODLR matrix $Y_{A,11}$ and a factorized low-rank matrix $Y_{A,21}$. We then update the second block column of $\tilde H$:
\begin{equation} \label{eq:update2ndblockH}
 \begin{bmatrix}
 \hat A_{12} \\
 \hat A_{22} \\
 \hat B_{R,2} \\
  \hat C_2
\end{bmatrix} := Q_1^T \begin{bmatrix}
 \tilde A_{12} \\
 \tilde A_{22} \\
 B_{R,2} \\
  C_2 
\end{bmatrix} = \begin{bmatrix}
 \tilde A_{12} - Y_{A,11} S \\
 \tilde A_{22} - Y_{A,21} S \\
  B_{R,2} - Y_{B_R,1} S \\ 
 C_2 - Y_{C,1} S \\
\end{bmatrix},
\end{equation}
where
\[
 S:= T^T_1 Y_1^T \begin{bmatrix}
 \tilde A_{12} \\
 \tilde A_{22} \\
  B_{R,2} \\
  C_2
\end{bmatrix} = T^T_1 ( Y_{A,11}^T \tilde A_{12} + Y_{A,21}^T \tilde A_{22} + Y_{B_R,1}^T B_{R,2} + Y_{C,1}^T C_2 ).
\]
It is important to note that each term of the sum in the latter expression is a low-rank matrix and, in turn, $S$ has low rank. This not only makes the computation of $S$ efficient but it also implies that the updates in~\eqref{eq:update2ndblockH} are of 
low rank and thus preserve the structure of the second block column of $\tilde H$.

After the update~\eqref{eq:update2ndblockH} has been performed, the process is completed by applying recursion to the updated second block column~\eqref{eq:update2ndblockH}, without the first block, and obtain a QR decomposition
\begin{equation}
\label{eq:second_block_QR}
 \begin{bmatrix}
 \hat{A}_{22} \\
 \hat B_{R,2} \\
  \hat C_2
\end{bmatrix} = Q_2 \begin{bmatrix}
  R_2 \\
  0 \\ 
\end{bmatrix},\quad Q_2 = I - Y_2 T_2 Y_2^T, \quad Y_2 = \begin{bmatrix} Y_{A,22}  \\ Y_{B_R,2} \\ Y_{C,2} \end{bmatrix}.
\end{equation}
By the discussion in Section~\ref{sec:recursive_qr}, see in particular~\eqref{eq:recqrcombined}, combining the QR decompositions of the first and second block columns yields a QR decomposition of $\tilde H$:
\begin{equation} \label{eq:tildeQRH}
  \tilde H =
 \tilde Q \begin{bmatrix} R \\ 0 \end{bmatrix} , \quad \tilde Q = I - \tilde YT \tilde Y^T
\end{equation}
with
\[
 \tilde Y = \left[ \begin{array}{c}
Y_{A} \\
\tilde Y_{B}  \\
Y_{C} 
\end{array} \right], \ 
Y_A = \left[ \begin{array}{c|c}
Y_{A,11} & 0 \\
Y_{A,21} & Y_{A,22}
\end{array} \right], \ 
\tilde Y_B = \left[ \begin{array}{c|c}
Y_{B_R,1} & Y_{B_R,2}
\end{array} \right], \ 
Y_C = \left[ \begin{array}{c|c}
Y_{C,1} & Y_{C,2}
\end{array} \right]
\]
and
\[
R = \left[ \begin{array}{c|c}
                      R_1 & \hat A_{12} \\  \hline
                      0 & R_2 
                      \end{array} \right], \quad T = \left[ \begin{array}{c|c}
     T_1 & -T_1 (Y_{A,21}^T Y_{A,22} + Y_{B_R,1}^T Y_{B_R,2} + Y_{C,1}^T Y_{C,2}) T_2 \\ \hline
     0 & T_2
\end{array} \right].
\]
Note that $Y_{A}$, $R$, and $T$ are triangular level-$\tilde \ell$ HODLR matrices, as desired.
\end{paragraph}

\begin{paragraph}{Postprocessing.} The procedure is completed by turning the obtained QR decomposition of $\tilde H$ into a QR decomposition of the matrix $H$ from~\eqref{eq:defH}. For 
this purpose, we simply set 
$Y_B = B_L \tilde Y_B$ and define $Y$ as in~\eqref{eq:defY}. Setting $Q = I-Y T Y^T$ then yields
\begin{eqnarray*}
  Q^T H &=& H - Y T^T Y^T H = \text{diag}(I_m,B_L,I_{r_2})\big( \tilde H - \tilde Y T^T \tilde Y^T \tilde H \big) \\
 &=&  \text{diag}(I_m,B_L,I_{r_2}) \tilde Q^T \tilde H = \text{diag}(I_m,B_L,I_{r_2}) \begin{bmatrix} R \\ 0 \end{bmatrix} = \begin{bmatrix} R \\ 0 \end{bmatrix}.
\end{eqnarray*}
Thus, we have obtained a QR decomposition of the form~\eqref{eq:qrdecomp}, which concludes the recursion step.
\end{paragraph}

\subsection{Algorithm and complexity estimates}
\label{sec:complexity_estimates}

Algorithm~\ref{alg:hqr} summarizes the recursive procedure described above. A QR decomposition of a level-$\ell$ HODLR matrix $A\in \R^{n\times n}$ is obtained by applying this 
algorithm with $\tilde A = A$ and void $B_L,B_R,C$.  

In the following, we derive complexity estimates for Algorithm~\ref{alg:hqr} applied to $A$ under the assumptions stated in Section~\ref{sec:HODLR}. In particular, it 
is assumed that all off-diagonal ranks (which are chosen adaptively) are bounded by $k$.

\begin{algorithm}[h!]
    \caption{\text{Recursive Householder based QR decomposition for HODLR matrices (hQR)}}
    \label{alg:hqr}
    \textbf{Input:} Level-$\tilde \ell$ HODLR matrix $\tilde{A}$, matrices $B_L,B_R,C$ defining 
    the matrix $H$ in~\eqref{eq:defH}. \\
    \textbf{Output:} Matrix $Y$ of the form~\eqref{eq:defY}, upper triangular 
    level-$\tilde l$ HODLR matrices $T,R$ defining an (approximate)  
    QR decomposition of $H = (I - YTY^T)\begin{bmatrix}R \\ 0 \end{bmatrix}$.   
    
\begin{algorithmic}[1]
\STATE {\textbf{function}} $[Y, T, R] =  {\tt hQR}(\tilde{A}, B, C)$
\IF {$B_L$ is not orthonormal} 
\STATE Compute economy-sized QR decomposition $B_L = Q R$ and set $B_L\gets Q, B_R \gets R B_R$. \label{left_orthogonal}
\ENDIF
 \STATE Set $\tilde{H} = \left[ \begin{smallmatrix} \tilde{A} \\ B_R \\  C\end{smallmatrix} \right]$. 
 \IF  {$\tilde \ell = 0$} 
\STATE  Use Alg.~\ref{alg:block_denseQR} to compute QR decomposition $\tilde{H} = (I - YTY^T)\begin{bmatrix}R \\ 0 \end{bmatrix}$ and partition $\tilde{Y} =\left[ \begin{smallmatrix} Y_A \\ \tilde{Y}_B \\ Y_C \end{smallmatrix} \right]$. \label{compute_lowest_level}
\ELSE   
\STATE Repartition $\tilde{H} =  \left[ \begin{smallmatrix}
 \tilde A_{11} & \tilde A_{12} \\
 \tilde A_{21} & \tilde A_{22} \\
  B_{R,1} &  B_{R,2} \\
  C_1 & C_2
\end{smallmatrix} \right]$ according to the HODLR format of $\tilde{A}$.
\STATE Compute QR decomposition of first block column of $\tilde{H}$ recursively: \\ $[Y_1, T_1, R_1] =  {\tt hQR} 
(\tilde{A}_{11}, \tilde{A}_{21}, \Big[ \begin{smallmatrix}  B_{R,1} \\ C_1 \end{smallmatrix} \Big] )$, 
with $Y_1$ defined by $Y_{A,11}$, $Y_{21}$,  $\Big[ \begin{smallmatrix} Y_{B_R,1} \\  Y_{C,1} \end{smallmatrix} \Big]$; see~\eqref{eq:first_block_QR}. \label{computeQR_first_blcol}         
\STATE Compute $\tilde{S} = \mathcal{T}_{\epsilon\cdot\Vert A\Vert_2} (Y_{A,11}^T \tilde A_{12} + Y_{A,21}^T \tilde A_{22} + Y_{B_R,1}^T B_{R,2} + Y_{C,1}^T C_2)$. \label{compute_Stilde} 
\STATE Compute $S =  T_1^T\tilde{S}$.\label{compute_S}
\STATE Update second block column of $\tilde{H}$:
$\left[ \begin{smallmatrix}
 \hat A_{12} \\
 \hat A_{22} \\
 \hat B_{R,2} \\
  \hat C_2
\end{smallmatrix} \right] := \left[ \begin{smallmatrix}
 \mathcal{T}_{\epsilon\cdot\Vert A\Vert_2} (\tilde A_{12} - Y_{A,11} S) \\
 \tilde A_{22} -_{\h} Y_{A,21} S \\
  B_{R,2} - Y_{B_R,1} S \\ 
 C_2 - Y_{C,1} S \\
\end{smallmatrix} \right] \text{.}$  \label{compute_update}
\STATE Compute QR decomposition of unreduced part of second block column of $\tilde{H}$ 
recursively: \\ $[Y_2, T_2, R_2] = {\tt hQR} (\hat A_{22},[\hskip 3pt], \Big[ \begin{smallmatrix} B_{R,2} \\ \hat C_2 \end{smallmatrix} \Big])$, with $Y_2$ defined by $Y_{A,22}$, $\Big[ \begin{smallmatrix} Y_{B_R,2} \\ Y_{C,2} \end{smallmatrix} \Big]$; see~\eqref{eq:second_block_QR}. \label{computeQR_second_blcol}
 \STATE Compute $\tilde T_{12} =\mathcal{T}_{\epsilon} ( Y_{A,21}^T Y_{A,22} + Y_{B_R,1}^T Y_{B_R,2} + Y_{C,1}^T Y_{C,2})$.  \label{compute_T}
\STATE Compute $T_{12} = -T_1 \tilde{T}_{12} T_2$. \label{compute_T2}
\STATE Set $T = \begin{bmatrix} T_{1} &T_{12} \\ 0 &T_{2}\end{bmatrix}$, and  $R = \begin{bmatrix} R_{1} &\hat{A}_{12} \\ 0 &R_{2} \end{bmatrix}$.
\STATE Set $ Y_A = \left[ \begin{array}{c|c}
Y_{A,11} & 0 \\
Y_{A,21} & Y_{A,22}
\end{array} \right], \ 
\tilde Y_B = \left[ \begin{array}{c|c}
Y_{B_R,1} & Y_{B_R,2}
\end{array} \right], \text{ and }\ 
Y_C = \left[ \begin{array}{c|c}
Y_{C,1} & Y_{C,2}
\end{array} \right].$
\ENDIF 
\STATE  Return $T, R$ and $Y= \left[ \begin{smallmatrix} Y_A \\ B_L \tilde{Y}_B \\ Y_C \end{smallmatrix} \right]$. 
\STATE {\textbf {end function}}
\end{algorithmic}
\end{algorithm} 

\paragraph{Line~\ref{left_orthogonal}. }
Every lower off-diagonal block of $A$ needs to be transformed once to left-orthogonal form in the course of the algorithm. For each $\tilde \ell$, $1\le \tilde \ell \le \ell$, there are $2^{\ell-\tilde \ell}$ such blocks of size $\calO(2^{\tilde \ell-1})\times \calO(2^{\tilde \ell-1})$ and rank at most $k$. Using the procedure for left-orthogonalization 
explained in Section~\ref{sec:hodlr}, the overall cost is
\[
\sum_{\tilde{\ell} = 1}^{\ell} \calO\big(2^{\ell-\tilde \ell} 2^{\tilde{\ell}-1} k^2\big) = \calO(k^2 n\log n),  
\]
where we used $\ell = \calO(\log n)$.

\paragraph{Line~\ref{compute_lowest_level}.} The QR decomposition of the compressed block column is performed for all $2^{\ell}$ block columns on the lowest level of recursion. Each of them is of size $\calO( \ell k) \times \calO(1)$, because there are at most $\ell$ lower off-diagonal blocks intersecting with each block column. Each QR decomposition requires $\calO( \ell k)$ operations and thus the overall cost is $\calO(k n\log n)$.

\paragraph{Lines~\ref{compute_Stilde}--\ref{compute_S}.} The computation of $S$ involves the following operations:
\begin{enumerate}
 \item three products of level-$(\tilde{\ell}-1)$ HODLR matrices with low-rank matrices;
 \item addition of four low-rank matrices, given in terms of their low-rank factors,  combined with recompression. 
\end{enumerate}
The first part is effected by performing at most $3k$ matrix-vector multiplications with $\calO(2^{\tilde \ell-1})\times \calO(2^{\tilde \ell-1})$ HODLR matrices. As the computation of $S$ is performed $2^{\ell - \tilde \ell}$ times for every $\tilde \ell$, we arrive at a total cost of
\begin{equation}
\label{eq:compute_Stilde1}
\sum_{\tilde{\ell}= 1}^{\ell} \calO\big( 2^{\ell - \tilde \ell} k^2 2^{\tilde \ell - 1} \tilde \ell \big) = \calO(k^2 n\log^2 n).
\end{equation}
In the second part, the addition is performed for matrices of size $\calO(2^{\tilde \ell-1})\times \calO(2^{\tilde \ell-1})$. The first three terms in Line~\ref{compute_Stilde} 
have rank at most $k$ but the rank of the last term $Y_{C,1}^T C_{2}$ can be up to $(\ell-1)k$. Letting the rank 
grow to $\calO(\ell k)$ would lead to an unfavourable complexity, because the cost of recompression depends quadratically on the rank of the matrix 
to be recompressed. To avoid this effect, we execute $(\ell + 1)k$ separate additions, each immediately followed by the
application of $\calT_{\epsilon\cdot \Vert A\Vert_2}$. Assuming that each recompression truncates to rank $\calO(k)$, this requires $\calO(2^{\tilde \ell-1} \ell k^2 )$ operations. Similarly as in~\eqref{eq:compute_Stilde1}, this leads to a total cost of $\calO(k^2 n\log^2 n)$.

\paragraph{Line~\ref{compute_update}.} For updating the second block column of $\tilde{H}$, the following operations are performed:
 \begin{enumerate}
 \item The computation of $\hat{A}_{12}$ requires $k$ HODLR matrix-vector multiplications, followed by low-rank recompression of a matrix of rank at most $2k$. Analogously to~\eqref{eq:compute_Stilde1}, this requires a total cost of $\calO(k^2 n \log^2n)$.

 \item The computation of $\hat{A}_{22}$ requires (approximate) subtraction of the product of two low-rank matrices from a level-$(\tilde \ell -1)$ HODLR matrix. The most expensive part of this step is the recompression of the updated HODLR matrix with off-diagonal ranks at most $2k$, amounting to a total cost of
\[\sum_{\tilde{\ell}= 1}^{\ell} \calO\big(2^{\ell-\tilde \ell} k^2 2^{\tilde \ell-1} \log 2^{\tilde \ell-1} \big) = \calO(k^2 n\log^2 n).
\]
\item The computation of $\hat{B}_{R,2}$ and $\hat{C}_{R,2}$ involves the multiplication of a matrix with at most $k+k\ell$ rows with a low-rank matrix, which requires $\calO(k^2 \ell 2^{\tilde \ell -1})$ operations each time. The total cost is thus again $\calO(k^2 n \log^2n)$.

\end{enumerate}

\paragraph{Lines~\ref{compute_T}--\ref{compute_T2}.} The computation of the low-rank block $T_{12}$ involves: 
\begin{enumerate}
 \item three multiplications of level-$(\tilde{\ell}-1)$ HODLR matrices with low-rank matrices;
 \item addition of three low-rank matrices, given in terms of their low-rank factors, combined with recompression. 
\end{enumerate}

Therefore, total cost is identical with the cost for computing $S$: $\calO(k^2 n \log^2n)$.

\paragraph{Summary.} The total cost of Algorithm~\ref{alg:hqr} applied to an $n\times n$ HODLR matrix is $\calO(k^2n \log^2 n)$.

\section{Numerical results} \label{sec:numericalresults}

In this section we demonstrate the efficiency of our method on several examples.  All algorithms 
were implemented and executed in \Matlab{} version R2016b on a dual Intel Core i7-5600U 2.60GHz CPU, $256$ KByte of level 2 cache and $12$ GByte of RAM, using a 
single core. Because all algorithms are rich in calls to BLAS and LAPACK routines, we believe that 
the use of \Matlab{} (instead of a compiled language) does not severely limit the predictive value of the reported timings.

The following algorithms have been compared:
\begin{description}
 \item[\bf CholQR] The Cholesky-based QR decomposition for HODLR matrices explained in Section~\ref{sec:cholesky}.
 \item[\bf CholQR2] CholQR followed by one step of the reorthogonalization procedure explained in Section~\ref{sec:cholesky}. 
 \item[\bf hQR] Algorithm~\ref{alg:hqr}, our newly proposed algorithm.
 \item[\bf MATLAB] Call to the MATLAB function {\tt qr}, which in turn calls the corresponding LAPACK routine for computing the QR decomposition of a general dense matrix.
\end{description}
\noindent Among the methods discussed in Section~\ref{sec:overview}, we have decided to focus on CholQR and CholQR2, primarily because they are relatively straightforward to implement. A comparison of 
CholQR and CholQR2 with the other methods from Section~\ref{sec:overview} can be found in~\cite{BennerMach2010}.

If not stated otherwise, when working with HODLR matrices we have chosen the minimal block size $n_{\min} = 250$ and the truncation tolerance $\epsilon = 10^{-10}$.

To assess accuracy, we have measured the numerical orthogonality of $Q= I - YTY^T$ and the residual of the computed QR decomposition:
\begin{equation} \label{eq:accuracy}
 e_{\mathrm{orth}} = \Vert Q^TQ - I\Vert_2, \quad   e_{\mathrm{acc}} = \Vert QR - A\Vert_2.
\end{equation}

\begin{example}[\bfseries{Performance for random HODLR matrices}]
\label{ex:time_vs_n_random} 
\rm We first investigate the performance of our method for HODLR matrices of varying size constructed as follows. The diagonal blocks are random dense matrices and 
each off-diagonal block is a rank-one matrix chosen as the outer product of two random vectors.  From Figure~\ref{fig:time_vs_n_random}, one observes that the computational 
time of the hQR algorithm nicely matches the $\calO(n\log^2n)$ reference line, the complexity claimed in Section~\ref{sec:complexity_estimates}. 
Compared to the much simpler and as we shall see, less accurate CholQR method, our new method is approximately 
only two times slower, while CholQR2 is slower than hQR for $n \geq 10\,000$. Note that for $n \geq 256\,000$, the Cholesky factorization of $A^T A$ fails to 
complete due to lack of (numerical) positive definiteness and, in turn, both CholQR and  CholQR2 return with an error.

\begin{figure}[ht!]
 \begin{center}
 \includegraphics[width=0.49\textwidth]{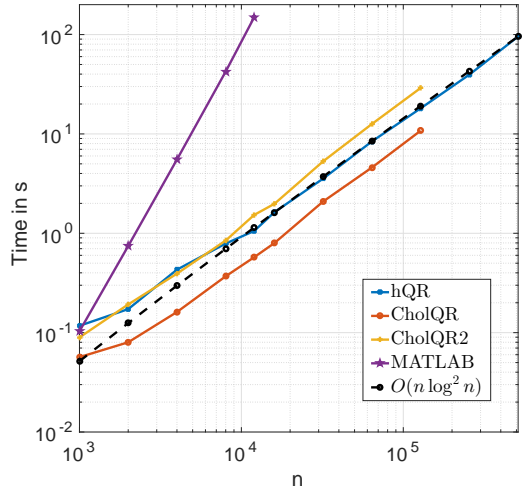}
 \end{center}
 \caption{Example~\ref{ex:time_vs_n_random}: Execution time vs. $n$ for computing QR decomposition of randomly generated $n\times n $ HODLR matrices. }
 \label{fig:time_vs_n_random}
 \end{figure}
 
Table~\ref{table:condition_vs_accuracy} provides insights into the observed accuracy for values of $n$ for which~\eqref{eq:accuracy} can be evaluated conveniently. As $n$ increases, the condition 
number of $A$ increases. Our method is robust to this increase and produces numerical orthogonality and a residual norm on the level of the truncation error. In contrast, the accuracy of CholQR clearly deteriorates as the condition number increases and the refinement performed by the more expensive CholQR2 cannot fully make up for this.

\begin{table}[ht!]
\caption{Example~\ref{ex:time_vs_n_random}. Numerical orthogonality $e_{\mathrm{orth}}$ and residual norm $e_{\mathrm{acc}}$, see~\eqref{eq:accuracy}, of different methods for computing QR decomposition of randomly generated $n\times n$ HODLR matrices.}
\label{table:condition_vs_accuracy}
\centering
{\renewcommand{\arraystretch}{1.2}
\begin{tabular}{c||c||c|c|c||c|c|c}
 $n$  &$\kappa_2(A)$  &$e^{\mathrm {hQR}}_{\mathrm{orth}}$ &$e^{\mathrm{CholQR}}_{\mathrm{orth}}$ &$e^{\mathrm{CholQR2}}_{\mathrm{orth}}$ &$e^{\mathrm{hQR}}_{\mathrm{acc}}$  &$e^{\mathrm{CholQR}}_{\mathrm{acc}}$   &$e^{\mathrm{CholQR2}}_{\mathrm{acc}}$ \\
 \hline
$1\,000$ &$8.7\cdot 10^4$  &$7.5\cdot 10^{-15}$  &$2.2\cdot 10^{-9}$ &$1.2\cdot 10^{-10} $ &$8.3\cdot 10^{-13}$  &$3.8\cdot 10^{-13}$ &$7.1\cdot 10^{-11}$  \\
\hline
$2\,000$ &$1.4\cdot 10^5$ &$1.4\cdot 10^{-14} $  &$3.4\cdot 10^{-8}$  &$2.4\cdot 10^{-9}$ &$4.4\cdot 10^{-12}$ &$1.1\cdot 10^{-12}$ &$1.3\cdot 10^{-9}$\\
 \hline
 $4\,000$ &$1.2\cdot 10^6$ &$1.6\cdot 10^{-13} $  &$8.4\cdot 10^{-7}$ &$1.5\cdot 10^{-8}$ &$1.5\cdot 10^{-11}$ &$7.6\cdot 10^{-12}$ &$1.6\cdot 10^{-8}$\\
 \hline
 $8\,000$ &$3.1\cdot 10^7$ &$1.9\cdot 10^{-12}$  &$5.1\cdot 10^{-6}$  &$5.7\cdot 10^{-8}$ &$1.9\cdot 10^{-10}$ &$9.1\cdot 10^{-10}$  &$3.5\cdot 10^{-8}$\\
 \hline 
$12\,000$  &$1.5\cdot 10^{8}$ &$1.8\cdot 10^{-12}$  &$2.2\cdot 10^{-4}$ &$8.2\cdot 10^{-7}$ &$1.9\cdot 10^{-10}$  &$1.4\cdot 10^{-10}$ &$2.7\cdot 10^{-7}$\\
\hline
 \end{tabular}
 }
\end{table}

Table~\ref{table:ranks_QR} aims at clarifying whether the representation of $Q$ in terms of its compact WY representation constitutes a disadvantage in terms of HODLR ranks. It turns out that the contrary 
is true; the maximal off-diagonal ranks of $Y$ and $T$ are significantly smaller than those of $Q$. Note, however, that does not translate into reduced memory consumption for the matrix sizes under 
consideration, because the larger off-diagonal ranks only occur in a few (smaller) off-diagonal blocks in $Q$.

\begin{table}[ht!]
\caption{Example~\ref{ex:time_vs_n_random}. Maximal off-diagonal ranks for the factors $Y, T, R$ and $Q = I - YTY^T$ from the QR decomposition computed by hQR 
applied to randomly generated $n\times n$ HODLR matrices. Memory for storing $Y$ and $T$ as well as $Q$ relative to memory for storing $A$ in the HODLR format. }
\label{table:ranks_QR}
\centering
{\renewcommand{\arraystretch}{1.1}
\begin{tabular}{c||c|c|c|c||c|c||}
 & \multicolumn{4}{c||}{Maximal ranks} & \multicolumn{2}{c||}{Memory} \\
 $n$  & $Y$ & $T$ & $Q$ & $R$ & $Y$ and $T$ & $Q$  \\
 \hline
$1\,000$ &$2$  &$2$ &$6$  &$4$ &$1.99$ &$1$ \\
\hline
$8\,000$ &$5$ &$5$  &$18$ &$10$ &$2$ &$1.2$  \\
 \hline
 $64\,000$ &$8$ &$8$ &$30$ &$15$ &$2.1$ &$1.5$ \\
 \hline
 $256\,000$ &$10$  &$10$ &$38$ &$17$ &$2.17$ &$1.7$ \\
 \hline
 \end{tabular}
 }
\end{table}

We also note that the maximal off-diagonal ranks and relative memory for $Y, T,R$ grow slowly, possibly logarithmically, as $n$ increases.

\end{example} 

\begin{example}[\bfseries{Accuracy for Cauchy matrices}]
\label{ex:accuracy_Cauchy} 
\rm In this example, we consider Cauchy matrices of size $n = 2000$, for which the entry $(i,j)$ is given by
$(x_i - y_j)^{-1}$ for $x,y \in \R^n$.
The vectors $x$ and $y$ are chosen as $2000$ equally spaced points from intervals $I_x$ and $I_y$, respectively, additionally perturbed by $\pm 2\cdot 10^{-2}$ with the sign chosen at random. We have used the following configurations:
\begin{itemize} 
 \item matrix $A_1$: intervals $I_x = [-1.25, 998.25]$ and $I_y = [-0.7, 998.9]$;  
 \item matrix $A_2$: intervals $I_x = [-1.25, 998.25]$ and $I_y = [-0.45, 999.15]$;  
 \item matrix $A_3$: intervals $I_x = [-1.25, 998.25]$ and $I_y = [-0.15, 999.45]$. 
\end{itemize}
All three matrices are invertible but their condition numbers are different. For each $i$, the HODLR approximation of $A_{i}$ has maximal off-diagonal rank $20$. Table~\ref{table:accuracy_Cauchy} summarizes 
the obtained results, which show that our method consistently attains an accuracy up to the level of truncation error. Once again, CholQR and CholQR2 fail to complete the computation for $A_3$, the most 
ill-conditioned matrix. In contrast to Example~\ref{ex:time_vs_n_random}, the maximal off-diagonal ranks for $Y$ and $T$ do not grow; they are bounded by $20$. 
The maximal off-diagonal rank for $R$ is $32$.
\begin{table}[ht!]
\caption{Example~\ref{ex:accuracy_Cauchy}. Numerical orthogonality $e_{\mathrm{orth}}$ and residual norm $e_{\mathrm{acc}}$, see~\eqref{eq:accuracy}, of different methods for 
computing QR decomposition of Cauchy matrices with varying condition number.}
\label{table:accuracy_Cauchy}
\centering
{\renewcommand{\arraystretch}{1.1}
\begin{tabular}{c||c||c|c|c||c|c|c}
  &$\kappa_2(A_i)$  &$e^{\mathrm {hQR}}_{\mathrm{orth}}$ &$e^{\mathrm{CholQR}}_{\mathrm{orth}}$ &$e^{\mathrm{CholQR2}}_{\mathrm{orth}}$ &$e^{\mathrm{hQR}}_{\mathrm{acc}}$  &$e^{\mathrm{CholQR}}_{\mathrm{acc}}$   &$e^{\mathrm{CholQR2}}_{\mathrm{acc}}$ \\
 \hline
 $A_1$   &$4.8\cdot 10^5$ &$5.7\cdot 10^{-11}$ &$2.6\cdot 10^{-5}$ &$2.8\cdot 10^{-11}$ &$1.1\cdot 10^{-8}$  &$2.9\cdot 10^{-9}$ &$1.9\cdot 10^{-7}$  \\
\hline
$A_2$  &$1.3\cdot 10^8$  &$3.6\cdot 10^{-10} $ &$1.3\cdot 10^{-1}$ &$3.4\cdot 10^{-9}$  &$2.3\cdot 10^{-9}$ &$6.8\cdot 10^{-10}$ &$4.6\cdot 10^{-9}$\\
 \hline
  $A_3$   &$2.9\cdot 10^{12}$  &$1.5\cdot 10^{-10} $ &-  &- &$2.2\cdot 10^{-9}$ &- &- \\
 \end{tabular}
 }
\end{table}
\end{example}

\begin{example}[\bfseries{Accuracy and orthogonality versus truncation tolerance}]
\label{ex:performance_vs_tolerance}
\rm 
As our final example, we investigate the influence of the truncation tolerance on the accuracy attained by our method. For this purpose, we consider the matrix $A_3$ from Example~\ref{ex:accuracy_Cauchy}. The 
truncation tolerance $\epsilon$  
is varied from $10^{-2}$ to $10^{-20}$, and the obtained results are compared with {\sc MATLAB}'s built-in function {\tt qr}. Figure~\ref{fig:performance_vs_tolerance} demonstrates that the 
errors decrease nearly proportional with $\epsilon$ until they stagnate around $\epsilon = 10^{-14}$, below which roundoff error appears to dominate.   

\begin{figure}[ht!]
 \begin{center}
 \includegraphics[width=0.49\textwidth]{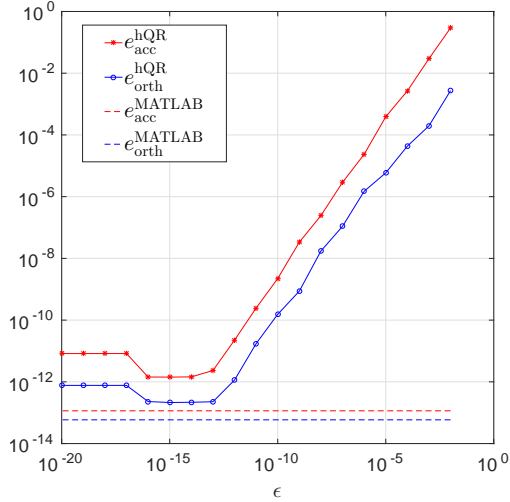}
 \end{center}
 \caption{Example~\ref{ex:performance_vs_tolerance}. Numerical orthogonality $e_{\mathrm{orth}}$ and residual norm $e_{\mathrm{acc}}$, see~\eqref{eq:accuracy}, of hQR applied to a Cauchy matrix with condition number $\approx 10^{12}$ vs. truncation tolerance $\epsilon$.}
 \label{fig:performance_vs_tolerance}
 \end{figure}

\end{example}

\section{Extension to rectangular HODLR matrices} \label{sec:rectangular}

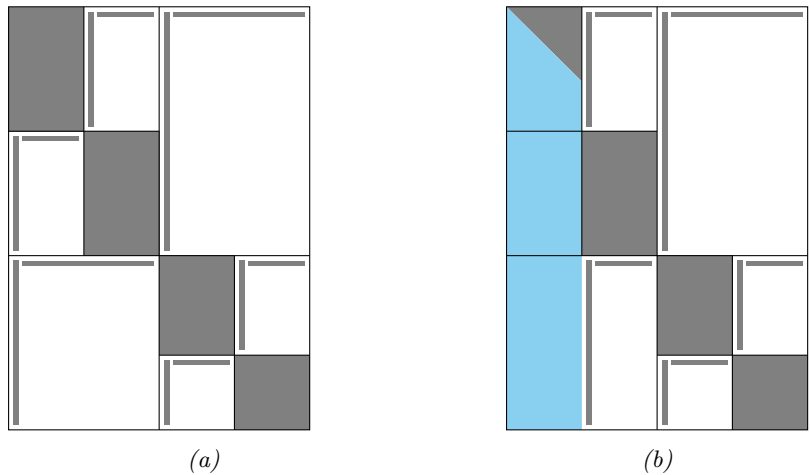
\begin{figure}[ht!]
 \centering
\begin{subfigure}[b]{0.35\textwidth}

\begin{tikzpicture}[scale = 0.33]
\fill[gray] (0,12)--(3,12)--(3,17)--(0,17);
\fill[gray] (3,7)--(6,7)--(6,12)--(3,12);

\fill[gray] (6,3)--(9,3)--(9,7)--(6,7);
\fill[gray] (9,0)--(12,0)--(12,3)--(9,3);

\fill [gray] (0.2, 0.2)--(0.4, 0.2)--(0.4,6.8)--(0.2, 6.8);
\fill [gray] (0.55, 6.8)--(0.55, 6.6)--(5.8, 6.6)--(5.8, 6.8);

\fill [gray] (6.2, 7.2)--(6.4, 7.2)--(6.4,16.8)--(6.2, 16.8);
\fill [gray] (6.55, 16.8)--(6.55, 16.6)--(11.8, 16.6)--(11.8, 16.8);

\fill [gray] (0.2, 7.2)--(0.4, 7.2)--(0.4,11.8)--(0.2, 11.8);
\fill [gray] (0.55, 11.8)--(0.55, 11.6)--(2.8, 11.6)--(2.8, 11.8);

\fill [gray] (3.2, 12.2)--(3.4, 12.2)--(3.4,16.8)--(3.2, 16.8);
\fill [gray] (3.55, 16.8)--(3.55, 16.6)--(5.8, 16.6)--(5.8, 16.8);

\fill [gray] (6.2, 0.2)--(6.4, 0.2)--(6.4,2.8)--(6.2, 2.8);
\fill [gray] (6.55, 2.8)--(6.55, 2.6)--(8.8, 2.6)--(8.8, 2.8);

\fill [gray] (9.2, 3.2)--(9.4, 3.2)--(9.4,6.8)--(9.2, 6.8);
\fill [gray] (9.55, 6.8)--(9.55, 6.6)--(11.8, 6.6)--(11.8, 6.8);

\draw (0,0) rectangle (12,17);  
\draw (0,7)--(12,7);
\draw (6,17)--(6,0);
\draw (0,12)--(6,12);
\draw (3,7)--(3,17);

\draw (6, 3)--(12,3);
\draw (9,0)--(9,7);
\end{tikzpicture}
\caption{}
\end{subfigure}
\qquad
\begin{subfigure}[b]{0.35\textwidth}
\centering
 \begin{tikzpicture}[scale = 0.33]
\fill[gray] (0,17)--(3,17)--(3,14);
\fill[gray] (3,7)--(6,7)--(6,12)--(3,12);

\fill[gray] (6,3)--(9,3)--(9,7)--(6,7);
\fill[gray] (9,0)--(12,0)--(12,3)--(9,3);

\fill [gray] (3.2, 0.2)--(3.4, 0.2)--(3.4,6.8)--(3.2, 6.8);
\fill [gray] (3.55, 6.8)--(3.55, 6.6)--(5.8, 6.6)--(5.8, 6.8);

\fill [gray] (6.2, 7.2)--(6.4, 7.2)--(6.4,16.8)--(6.2, 16.8);
\fill [gray] (6.55, 16.8)--(6.55, 16.6)--(11.8, 16.6)--(11.8, 16.8);


\fill [gray] (3.2, 12.2)--(3.4, 12.2)--(3.4,16.8)--(3.2, 16.8);
\fill [gray] (3.55, 16.8)--(3.55, 16.6)--(5.8, 16.6)--(5.8, 16.8);

\fill [gray] (6.2, 0.2)--(6.4, 0.2)--(6.4,2.8)--(6.2, 2.8);
\fill [gray] (6.55, 2.8)--(6.55, 2.6)--(8.8, 2.6)--(8.8, 2.8);

\fill [gray] (9.2, 3.2)--(9.4, 3.2)--(9.4,6.8)--(9.2, 6.8);
\fill [gray] (9.55, 6.8)--(9.55, 6.6)--(11.8, 6.6)--(11.8, 6.8);

\fill [color1] (0,0)--(3,0)--(3,14)--(0,17);

\draw (0,0) rectangle (12,17);  
\draw (0,7)--(12,7);
\draw (6,17)--(6,0);
\draw (0,12)--(6,12);
\draw (3,7)--(3,17);

\draw (6, 3)--(12,3);
\draw (9,0)--(9,7);
\end{tikzpicture}
\caption{}
\end{subfigure}
 \caption{ (a) Rectangular HODLR matrix of level $\ell = 2$. (b)
 Structure after reducing the first block column. Zero parts of the matrix are colored blue.}
 \label{fig:hodlr_rectangular}
\end{figure}

In this section, we sketch the extension of Algorithm~\ref{alg:hqr} to rectangular HODLR matrices. For such matrices, one allows the diagonal blocks in the 
recursive partitioning~\eqref{eq:toplevelpartitioning} to be rectangular. In turn, the definition of a  rectangular HODLR matrix $A\in \R^{m\times n}$ depends on two integer partitions
\[
m = m_1 + m_2 + \cdots m_{2^\ell}, \quad  n = n_1 + n_2 + \cdots n_{2^\ell},  
\]
corresponding to the sizes $m_j \times n_j$, $j = 1,\ldots, 2^\ell$, of the diagonal blocks on the lowest level of the recursion. In the following, we assume that
\[
 m_j \ge n_j, \quad j = 1,\ldots, 2^\ell.
\]
See Figure~\ref{fig:hodlr_rectangular} (a) for an illustration.

An application appears in our work~\cite{Susnjara2018} on a fast spectral divide-and-conquer method, which requires the QR decomposition of an $m\times n$ matrix that results from selecting $n\approx m/2$ columns of
an $m\times m$ HODLR matrix.

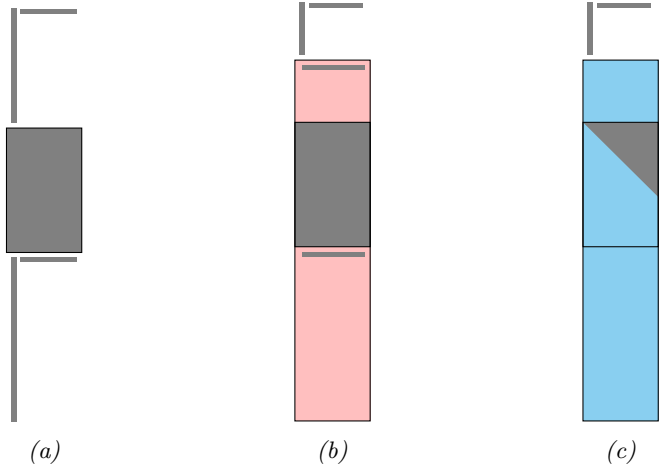
\begin{figure}[ht!]
 \centering
\begin{subfigure}[b]{0.25\textwidth}
\centering
 \begin{tikzpicture}[scale = 0.33]
\fill[gray] (3,7)--(6,7)--(6,12)--(3,12);

\fill [gray] (3.2, 0.2)--(3.4, 0.2)--(3.4,6.8)--(3.2, 6.8);
\fill [gray] (3.55, 6.8)--(3.55, 6.6)--(5.8, 6.6)--(5.8, 6.8);

\fill [gray] (3.2, 12.2)--(3.4, 12.2)--(3.4,16.8)--(3.2, 16.8);
\fill [gray] (3.55, 16.8)--(3.55, 16.6)--(5.8, 16.6)--(5.8, 16.8);

\draw (3,7) rectangle (6,12); 
\end{tikzpicture}
\caption{}
\end{subfigure}
\begin{subfigure}[b]{0.25\textwidth}
\centering
\begin{tikzpicture}[scale = 0.33]
\fill[pink] (3, 0)--(6,0)--(6,14.5)--(3, 14.5);
\fill[gray] (3,7)--(6,7)--(6,12)--(3,12);

\fill [gray] (3.3, 6.8)--(3.3, 6.6)--(5.8, 6.6)--(5.8, 6.8);
\fill [gray] (3.3, 14.3)--(3.3, 14.1)--(5.8, 14.1)--(5.8, 14.3);

\fill [gray] (3.2, 14.7)--(3.4, 14.7)--(3.4,16.8)--(3.2, 16.8);
\fill [gray] (3.55, 16.8)--(3.55, 16.6)--(5.7, 16.6)--(5.7, 16.8);

\draw (3,7) rectangle (6,12); 
\draw (3,0) rectangle (6,14.5); 

\end{tikzpicture}
\caption{}
\end{subfigure}
\begin{subfigure}[b]{0.25\textwidth}
\centering
\begin{tikzpicture}[scale = 0.33]
\fill[color1] (3,0)--(6,0)--(6,14.5)--(3,14.5);

\fill[gray] (3,12)--(6,12)--(6,9);

\fill [color1] (3.2, 0.2)--(3.4, 0.2)--(3.4,6.8)--(3.2, 6.8);
\fill [color1] (3.3, 6.8)--(3.3, 6.6)--(5.8, 6.6)--(5.8, 6.8);
\fill [color1] (3.3, 14.3)--(3.3, 14.1)--(5.8, 14.1)--(5.8, 14.3);

\fill [gray] (3.2, 14.7)--(3.4, 14.7)--(3.4,16.8)--(3.2, 16.8);
\fill [gray] (3.55, 16.8)--(3.55, 16.6)--(5.7, 16.6)--(5.7, 16.8);

\draw (3,7) rectangle (6,12); 
\draw (3,0) rectangle (6,14.5); 
\end{tikzpicture}
\caption{}
\end{subfigure}
 \caption{Reduction of the second block column of a rectangular HODLR matrix of level $\ell = 2$. Zero parts are colored blue.}
 \label{fig:hodlr_2ndcol}
\end{figure}

We now consider the application of Algorithm~\ref{alg:hqr} to a rectangular HODLR matrix. This algorithm starts with reducing the first block column to upper triangular form. This first step is coherent 
with the structure, see Figure~\ref{fig:hodlr_rectangular} (b), and no significant modification of Algorithm~\ref{alg:hqr} is necessary. However, the same cannot be said about the subsequent steps. 
The transformation of the second block column (or, more precisely, its unreduced part) to upper triangular form would mix dense with 
low-rank blocks and in turn destroy the HODLR format. To avoid this effect, we reduce the second block column to \emph{permuted triangular form}, such that the 
reduced triangular matrix replaces the dense diagonal block and all other parts become zero. This process is illustrated in Figure~\ref{fig:hodlr_2ndcol}: First the low-rank blocks in the unreduced part (the part highlighted in Figure~\ref{fig:hodlr_2ndcol} (b)) are compressed. Then an orthogonal transformation is performed such that the nonzero rows are reduced to a triangular matrix situated on top of the dense block; see Figure~\ref{fig:hodlr_2ndcol} (c). In practice, this is effected by an appropriate permutation of the rows, followed by a QR decomposition and the inverse permutation. The rows of the factor $Y$ in the compact WY representation of this transformation are permuted accordingly and, in turn, $Y$ inherits the structure from the second block column.

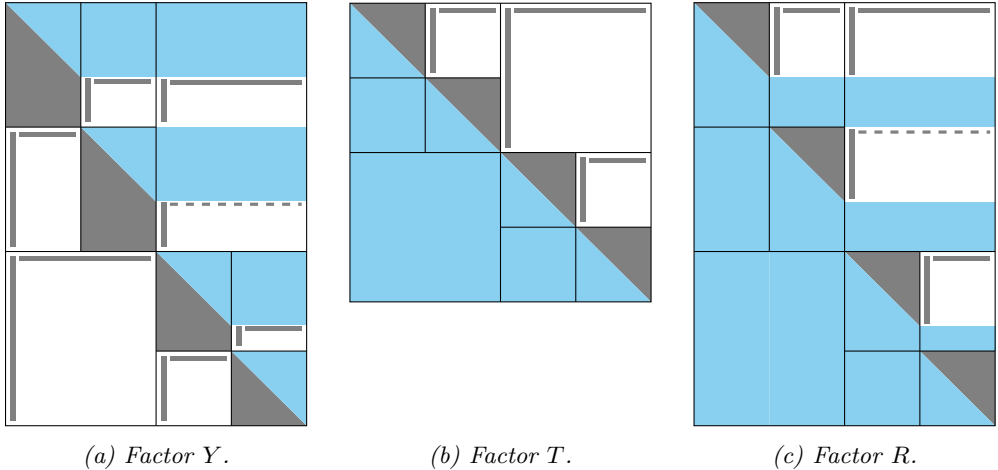
\begin{figure}[ht!]
 \centering

\begin{subfigure}[b]{0.3\textwidth}
\centering
\begin{tikzpicture}[scale = 0.33]
\fill[gray] (0,12)--(3,12)--(3,14)--(0,17);
\fill[gray] (3,7)--(6,7)--(6,9)--(3,12);

 \fill[gray] (6,3)--(9,3)--(9,4)--(6,7);
\fill[gray] (9,0)--(12,0)--(9,3);
\fill [gray] (0.2, 0.2)--(0.4, 0.2)--(0.4,6.8)--(0.2, 6.8); 
\fill [gray] (0.55, 6.8)--(0.55, 6.6)--(5.8, 6.6)--(5.8, 6.8); 

\fill [gray] (6.2, 7.2)--(6.4, 7.2)--(6.4,9)--(6.2, 9); 
\fill [gray] (6.2, 12.2)--(6.4, 12.2)--(6.4,13.9)--(6.2, 13.9); 
\fill [gray] (6.55, 13.9)--(6.55, 13.7)--(11.8, 13.7)--(11.8, 13.9); 

\fill [gray] (0.2, 7.2)--(0.4, 7.2)--(0.4,11.8)--(0.2, 11.8); 
\fill [gray] (0.55, 11.8)--(0.55, 11.6)--(2.8, 11.6)--(2.8, 11.8); 

\fill [gray] (3.2, 12.2)--(3.4, 12.2)--(3.4,13.95)--(3.2, 13.95); 
\fill [gray] (3.55, 13.95)--(3.55, 13.75)--(5.8, 13.75)--(5.8, 13.95); 

\fill [gray] (6.2, 0.2)--(6.4, 0.2)--(6.4,2.8)--(6.2, 2.8); 
\fill [gray] (6.55, 2.8)--(6.55, 2.6)--(8.8, 2.6)--(8.8, 2.8); 

\fill [gray] (9.2,3.2)--(9.4,3.2)--(9.4,4)--(9.2,4); 
\fill [gray] (9.55,4)--(9.55,3.8)--(11.8, 3.8)--(11.8,4); 


\fill [color1] (3,14.01)--(12,14.01)--(12,17)--(3,17);
\fill [color1] (6,9.05)--(12,9.05)--(12,12)--(6,12);
\fill [color1] (9,4.05)--(12,4.05)--(12,7)--(9,7);
\fill [color1] (0,17)--(3,17)--(3,14); 
\fill [color1] (3,12)--(6,12)--(6,9);
\fill [color1] (6,7)--(9,7)--(9,4);
\fill [color1] (9,3)--(12,3)--(12,0);

\draw (0,0) rectangle (12,17);  
\draw (0,7)--(12,7);
\draw (6,17)--(6,0);
\draw (0,12)--(6,12);
\draw (3,7)--(3,17);

\draw (6, 3)--(12,3);
\draw (9,0)--(9,7);

\draw [gray, dashed, very thick] (6.55, 8.9)--(11.8, 8.9); 
\end{tikzpicture}
\caption{Factor $Y$.}
\end{subfigure}
\begin{subfigure}[b]{0.3\textwidth}
\centering
\begin{tikzpicture}[scale = 0.33]

\fill [gray] (0,12)--(3,12)--(3,9);
\fill [gray] (3,9)--(6,9)--(6,6);
\fill [gray] (6,6)--(9,6)--(9,3);
\fill [gray] (9,3)--(12,3)--(12,0);

\fill [gray] (6.2,6.2)--(6.4, 6.2)--(6.4,11.8)--(6.2, 11.8);
\fill [gray] (6.55,11.8)--(6.55, 11.6)--(11.8,11.6)--(11.8, 11.8);

\fill [gray] (3.2,9.2)--(3.4, 9.2)--(3.4,11.8)--(3.2, 11.8);
\fill [gray] (3.55,11.8)--(3.55, 11.6)--(5.8,11.6)--(5.8, 11.8);

\fill [gray] (9.2,3.2)--(9.4, 3.2)--(9.4,5.8)--(9.2, 5.8);
\fill [gray] (9.55,5.8)--(9.55, 5.6)--(11.8,5.6)--(11.8, 5.8);

\fill [color1] (0,0)--(12,0)--(0,12);

\fill [white] (0,-5)--(3,-5)--(3,0)--(0,0);

\draw  (0,0) rectangle (12,12);
\draw  (6,0)--(6,12);
\draw (0,6)--(12,6);

\draw (3,6)--(3,12);
\draw (0,9)--(6,9);
\draw (9,0)--(9,6);
\draw (6,3)--(12,3);

\end{tikzpicture}
\caption{Factor $T$.}
\end{subfigure}
\begin{subfigure}[b]{0.3\textwidth}
\centering
 \begin{tikzpicture}[scale = 0.33]
\fill[gray] (0,17)--(3,17)--(3,14);
\fill[gray] (3,12)--(6,12)--(6,9);

\fill[gray] (6,7)--(9,7)--(9,4);
\fill[gray] (9,3)--(12,3)--(12,0);

\fill [gray] (6.2, 14.1)--(6.4, 14.1)--(6.4,16.8)--(6.2, 16.8);
\fill [gray] (6.2, 9.1)--(6.4, 9.1)--(6.4,11.9)--(6.2, 11.9);
\fill [gray] (6.55, 16.8)--(6.55, 16.6)--(11.8, 16.6)--(11.8, 16.8);

\fill [gray] (3.2, 14.1)--(3.4, 14.1)--(3.4,16.8)--(3.2, 16.8); 
\fill [gray] (3.55, 16.8)--(3.55, 16.6)--(5.8, 16.6)--(5.8, 16.8); 

\fill [gray] (9.2, 4.1)--(9.4, 4.1)--(9.4,6.8)--(9.2, 6.8);
\fill [gray] (9.55, 6.8)--(9.55, 6.6)--(11.8, 6.6)--(11.8, 6.8);

\fill [color1] (0,0)--(3,0)--(3,14)--(0,17);
\fill [color1] (3,0)--(6,0)--(6,9)--(3,12);
\fill [color1] (6,0)--(9,0)--(9,4)--(6,7);
\fill [color1] (9,0)--(12,0)--(9,3);

\fill [color1] (9,3)--(12,3)--(12,4)--(9,4);
\fill [color1] (3,12)--(12,12)--(12,14)--(3,14);
\fill [color1] (6,7)--(12,7)--(12,9)--(6,9);

\draw (0,0) rectangle (12,17);  
\draw (0,7)--(12,7);
\draw (6,17)--(6,0);
\draw (0,12)--(6,12);
\draw (3,7)--(3,17);

\draw (6, 3)--(12,3);
\draw (9,0)--(9,7);
\draw [gray, dashed, very thick] (6.55, 11.8)--(11.8, 11.8); 
\end{tikzpicture}
\caption{Factor $R$.}
\end{subfigure}
\caption{Illustration of factors $Y, T$ and $R$ of a (permuted) QR decomposition of a rectangular HODLR matrix of level $\ell = 2$. Dashed lines denote right low-rank factors that are shared with the
 off-diagonal blocks above the considered block.}
 \label{fig:rectQR}
\end{figure}

The described process is applied to each block column on the lowest level of the recursion: A permuted QR decomposition is performed such that 
the reduced $n_j\times n_j$ triangular matrix is situated on top of the dense diagonal block. On higher 
levels of the recursion, Algorithm~\ref{alg:hqr} extends with relatively minor modifications. This modified algorithm results in a 
QR decomposition $A = QR$, $Q = I -YTY^T$, where $Y,R$ are \emph{permuted} lower trapezoidal/upper triangular matrices that inherit the HODLR format 
of $A$. The matrix $T$ is an $n\times n$ upper triangular HODLR matrix; see Figure~\ref{fig:rectQR} for an illustration. Note that, in particular, $R$ is not triangular, but it can be easily permuted to triangular form, if needed.

We have collected preliminary numerical evidence that the described modified algorithm is effective at computing permuted QR decompositions of 
rectangular HODLR matrices. For this purpose, we have applied a dense version of the algorithm and compressed the obtained 
factors $Y$, $T$, $R$ afterwards, in accordance with the format shown in Figure~\ref{fig:rectQR}. The parameters guiding 
the HODLR format are identical to the default parameters in Section~\ref{sec:numericalresults}: $n_{\min} = 250$ and $\epsilon = 10^{-10}$.

\begin{example}[\bfseries{Performance for an invariant subspace basis}]
 \label{ex:hsdc_rectangular}
 \rm This example illustrates the use of our algorithm for orthonormalizing a set of vectors in an application from~\cite{Susnjara2018}.
For this purpose, we consider a tridiagonal symmetric matrix $T\in \R^{m\times m}$ with $m = 8\,000$, chosen such that the eigenvalues are uniformly 
distributed in $[-1, \hskip 2pt -10^{-1}]\cup [10^{-1}, \hskip 2pt 1]$. 
It turns out that the spectral projector $\Pi_{<0}$ associated with the negative eigenvalues of $T$ can be well approximated in the HODLR format; the numerical ranks of the 
off-diagonal blocks are bounded by $17$. We applied the 
method proposed in~\cite[Section 4.1]{Susnjara2018} (with threshold parameter $\delta = 0.35$) to select a 
well-conditioned set of $n\approx m/2$ columns of $\Pi_{<0}$, which will be denoted by $\Pi_{<0}(:, C) \in \R^{m\times n}$ with the column indices $C$.

We applied the described modification of Algorithm~\ref{alg:hqr} to orthonormalize the rectangular HODLR matrix $A = \Pi_{<0}(:, C) \in \R^{m\times n}$; an operation 
 needed in~\cite{Susnjara2018}. The accuracy we obtained is at the level of truncation tolerance: $e_{\mathrm{orth}} = 5.8\cdot 10^{-12}$ and $e_{\mathrm{acc}} = 8.9\cdot 10^{-11}$. The algorithm is 
also efficient in terms of memory; see Table~\ref{table:rectangular}. In particular, the off-diagonal ranks of the factors $Y,T$ and $R$ do not grow compared 
 to $A$. In this example, and in contrast to the square examples reported in Section~\ref{sec:numericalresults}, the memory is reduced when storing $Y$ and $T$ instead of $Q$.

\begin{table}[ht!]
\caption{Examples~\ref{ex:hsdc_rectangular} and~\ref{ex:random_rectangular}. Maximal off-diagonal ranks for the factors $Y, T, R$ and $Q = I - YTY^T$.  Memory for storing $Y$ and $T$ as 
well as $Q$, and $R$ relative to memory for storing $A$ in the HODLR format. }
\label{table:rectangular}
\centering
{\renewcommand{\arraystretch}{1.2}
\begin{tabular}{l||c|c|c|c||c|c|c||}
 & \multicolumn{4}{c||}{Maximal ranks} & \multicolumn{3}{c||}{Memory} \\
&$Y$ & $T$ & $Q$ & $R$ & $Y$ and $T$  &$Q$ & $R$ \\
\hline
Example~\ref{ex:hsdc_rectangular} &$14$ &$14$ &$23$ &$11$ &$1.5$  &$2.1$ &$0.87$\\
\hline
\hline
Example~\ref{ex:random_rectangular} &$8$ &$12$ &$12$ &$8$ &$1.6$ &$2.2$ &$1.1$ \\
\hline
\end{tabular} }
\end{table}
  
\end{example}

\begin{example}[\bfseries{Performance for a random rectangular HODLR matrix}]
 \label{ex:random_rectangular}
 \rm In analogy to Example~\ref{ex:time_vs_n_random} we  generated a random $8\,000 \times 4\,000$ HODLR matrix $A$ with 
 off-diagonal ranks $1$. We obtained 
 $e_{\mathrm{orth}} = 2.8\cdot 10^{-13}$, $e_{\mathrm{acc}} = 1.4\cdot 10^{-11}$, and the off-diagonal ranks and the memory 
requirements shown in Table~\ref{table:rectangular}.
\end{example}

 \section{Conclusion}

We have presented the hQR method, a novel, fast and accurate method for computing the QR decomposition of a HODLR matrix.
Our numerical experiments indicate that hQR is the method of choice, unless one wants to sacrifice accuracy for a relatively small gain in computational time.
It remains to be seen whether the developments of this work extend to the broader class of hierarchical matrices.

\bibliographystyle{plain}
\bibliography{biblio}

\end{document}